\input amstex
\input amsppt.sty

\magnification1200
\vsize19cm


\def\LL{\leavevmode\setbox0=\hbox{L}\hbox to\wd0{\hss\char'40L}}

\def\rh{\rho}

\def\ph{\varphi}


\def\R{{\Bbb R}}

\def\C{{\Bbb C}}

\def\N{{\Bbb N}}
\def\Z{{\Bbb Z}}
\def\Q{{\Bbb Q}}

\def\today{\ifcase\month\or
 January\or February\or March\or April\or May\or June\or
 July\or August\or September\or October\or November\or December\fi
 \space\number\day, \number\year}

\def\tfrac#1#2{{\textstyle{#1\over#2}}}

\def\({\left(}
\def\){\right)}
\def\[{\left[}
\def\]{\right]}

\def\suml{\sum\limits}

\def\prodl{\prod\limits}

\def\3{\ss}

\def\HT{\operatorname{ht}}

\topmatter 
\title
Orthogonal polynomials associated with root systems
\endtitle 
\author I. G. Macdonald 
\endauthor 
\affil 
School of Mathematical Sciences, Queen Mary and Westfield College, \\
University of London, London E1 4NS, England
\endaffil
\address 
School of Mathematical Sciences, Queen Mary and Westfield College, 
University of London, London E1 4NS, England
\endaddress
\abstract Let $R$ and $S$ be two irreducible root systems spanning
the same vector space and having the same Weyl group $W$, such that
$S$ (but not necessarily $R$) is reduced. For each such pair $(R,S)$
we construct a family of $W$-invariant orthogonal polynomials in
several variables, whose coefficients are rational functions of
parameters $q,t_1,t_2,\dots,t_r$, where $r$ ($=1,2$ or 3) is the
number of $W$-orbits in $R$. For particular values of these
parameters, these polynomials give the values of zonal spherical
functions on real and $p$-adic symmetric spaces. Also when $R=S$ is
of type $A_n$, they conincide with the symmetric polynomials
described in I.~G.~Macdonald, {\it Symmetric Functions and Hall
Polynomials}, 2nd edition, Oxford University Press (1995),
Chapter~VI.
\endabstract
\endtopmatter

\document

\parskip.1cm plus 3pt minus -3pt

\head Foreword\endhead
The text that follows this Foreword is that of my 1987 preprint with
the above title. It is now in many ways a period piece, and I have
thought it best to reproduce it unchanged. I am grateful to Tom
Koornwinder and Christian Krattenthaler for arranging for its
publication in the S\'eminaire Lotharingien de Combinatoire.

I should add that the subject has advanced considerably in the
intervening years. In particular, the conjectures in \S12 below are
now theorems. For a sketch of these later developments the reader may
refer to my booklet {\it ``Symmetric functions and orthogonal
polynomials"}, University Lecture Series Vol.~12, American
Mathematical Society (1998), and the references to the literature
given there.

\bigskip
\noindent
November 2000
\bigskip

\head \bf Introduction\endhead

The orthogonal polynomials which are the subject of this paper are
Laurent polynomials in several variables. To be a little more
precise, they are elements of the group algebra $A$ of the weight
lattice $P$ of a root system $R$, invariant under the action of the
Weyl group of $R$, and they depend rationally on two parameters $q$ and
$t$\hbox{\ }\footnote"(*)"{\hbox{\ }This is a simplified description for the purposes of
this introduction.}. They are indexed by the dominant weights and they
are pairwise orthogonal with respect to a certain weight function
$\triangle$, to be defined later.

For particular values of the parameters $q$ and $t$, these
polynomials reduce to familiar objects:

\medskip
(i) when $q=t$ they are independent of $q$ and are the Weyl
characters for the root system $R$.

\smallskip
(ii) when $t=1$ they are again independent of $q$, and are the
elements of $A$ corresponding to the orbits of the Weyl group in the
weight lattice $P$.

\smallskip
(iii) when $q=0$ they are (up to a scalar factor) the polynomials that
give the values of zonal spherical functions on a semisimple $p$-adic
Lie group $G$ relative to a maximal compact subgroup $K$, such that
the restricted root system of $(G,K)$ is the dual root system $R^\lor$.
Here the value of the parameter $t$ is the reciprocal of the
cardinality of the residue field of the local field over which $G$ is
defined.

\smallskip
(iv) finally, when $q$ and $t$ both tend to $1$, in such a way
that $(t-1)/(q-1)$ tends to a definite limit $k$, then (for certain
values of $k$) our polynomials give the values of zonal spherical
functions on a real (compact or non-compact) symmetric space $G/K$
arising from finite-dimensional representations of $G$ that have a
$K$-fixed vector $\ne 0$. Here the root system $R$ is the restricted
root system of $G/K$, and the parameter $k$ is half the root
multiplicity (assumed for the purposes of this description to be the
same for all restricted roots).

\medskip
Thus these two-parameter families of orthogonal polynomials constitute
a sort of bridge between harmonic analysis on real symmetric spaces
and on their $p$-adic analogues. It is perhaps natural to ask, in view
of recent developments (quantum groups, etc.), whether there is a
group-like object depending  on two parameters $q$ and $t$ that lies
behind this theory; but on this question we have nothing to say. We
would only remark that such a hypothetical object would have to
partake of the properties of a $p$-adic Lie group when $q=0$, and of a
real Lie group in the limiting case $(q,t)\rightarrow(1,1)$ described
in (iv) above.

All this is in fact a simplified description of the theory. The
context in which we shall work is that of an ``admissible pair" $(R,S)$
of root systems: this means that $R$ and $S$ are finite root systems in
the same vector space, having the same Weyl group $W$ and such that
$S$ (but not necessarily $R$) is reduced. In this context we define
parameters $q_\alpha$ and $t_\alpha$ for each root $\alpha\in R$, such
that $q_\alpha = q_\beta$ and $t_\alpha = t_\beta$ if $\alpha$ and
$\beta$ are in the same $W$-orbit. This is described, and the
appropriate notation established, in Sections~1 and 
2. The weight
function $\triangle$ and the accompanying scalar product are defined
in Section~3. The main result of the paper is Theorem~(4.1), which
asserts the existence of a family of orthogonal polynomials
$P_\lambda$ associated with a given admissible pair $(R,S)$.

The proof of the theorem consists in constructing a suitable
self-adjoint linear operator $E$ with distinct eigenvalues; the
polynomials $P_\lambda$ are the eigenfunctions of $E$, suitably
normalized. In fact we need two constructions for such a linear
operator. The first of these is described in Section~5, and works
whenever the root system $S$ (assumed irreducible) possesses a
minuscule weight, that is to say provided that $S$ is not of type
$E_8$, $F_4$ or $G_2$. The second construction, described in 
Section~6, is based on the premise, familiar to experts in standard monomial
theory, that the next best thing to a minuscule weight is a 
quasi-minuscule weight, and produces an operator $E$ with the desired
properties in the cases not covered by the previous
construction.

In Sections~8--11 we consider the particular cases corresponding to
(i)-(iv) above. We also consider, in Section~9, the case where $R$ is
of rank 1. If $R$ is of type $A_1$, the polynomials $P_\lambda$ are
essentially the $q$-ultraspherical polynomials \cite{1}, whereas if $R$ is of
type $BC_1$ the $P_\lambda$ reduce to a particular case of the
orthogonal polynomials of Askey and Wilson \cite{2}. Also, if $R$ is of
type $A_n$ ($n\geq1$) the $P_\lambda$ are essentially the symmetric
functions that are the subject of Chapter VI of \cite{11}.

Finally, in Section~12 we put forward two conjectures relating to the
polynomials $P_\lambda$. They involve a common generalization of
Harish-Chandra's $c$-function and its $p$-adic  counterpart, and one of
the conjectures includes as a special case the constant term
conjectures of \cite{10} and \cite{13}.

\head \bf \S1 \endhead

Let $V$ be a real vector space of finite dimension, endowed with a
positive-definite symmetric bilinear form $\langle  u, v \rangle$. We
shall write $\vert v\vert = \langle v, v \rangle ^{1/2}$ for $v\in V$, and
$$v^\lor =2v/\vert v\vert^2$$
if $v\neq 0$. If $R$ is a root system in $V$, we denote by $R^\lor $ the
dual root system $\{\alpha^\lor :\alpha \in R\}$.

Let $R$ and $S$ be root systems in $V$ (and spanning $V$). The pair
$(R,S)$ will be said to be  {\it admissible} if $R$ and $S$ have
the same Weyl group $W$, and $S$ (but not necessarily $R$) is
reduced.

Suppose that $(R,S)$ is admissible. Then the set of hyperplanes in $V$
orthogonal to the roots is the same for both $R$ and $S$, and hence
(as $S$ is reduced) there exists for each $\alpha\in R$ a unique positive real
number $u_\alpha$ such that
$$\alpha_\ast=u^{-1}_\alpha \alpha\in S,$$
and the mapping $f:R\rightarrow S$ defined by $f(\alpha)=\alpha_\ast$
is surjective.

Let $\alpha\in R$, $w\in W$ and let $\beta = w\alpha$. Then
$w(\alpha_\ast)\in S$, and is a positive scalar multiple of $\beta$,
so that $w(\alpha_\ast)=\beta_\ast = (w\alpha)_\ast$. Hence the
mapping $f$ commutes with the action of $W$, and $u_\alpha = u_\beta$
whenever $\alpha$, $\beta$ lie in the same $W$-orbit in $R$. Moreover
if $R$ is not reduced and $\alpha$, $2\alpha\in R$, we have
$$u_{2\alpha} = 2u_\alpha,\tag 1.1$$ 
since $(2\alpha)_\ast = \alpha_\ast$.

From now on we shall assume that $R$ (and therefore also $S$) is
{\it irreducible}. Another pair $(R',S')$ of root systems in $V$
will be said to be  {\it similar} to $(R,S)$ if there exist
positive real numbers $a,b$ such that $R' = aR$ and $S' = bS$. The
effect of passing from $(R,S)$ to a similar pair is simply to multiply
each $u_\alpha$ by the same positive scalar factor.

The classification of irreducible admissible pairs $(R,S)$ up to
similarity is easily described. There are three cases to consider.

\medskip
(i) $R$ is reduced and $S=R$, so that $u_\alpha = 1$ for each
$\alpha\in R$.

\smallskip
(ii) $R$ is reduced, with two root-lengths, and $S=R^\lor $. Then
$u_\alpha = \frac {1} {2}\vert \alpha\vert^2$ for each $\alpha\in R$.
We may assume that $\vert\alpha\vert^2 = 2$ for each short root
$\alpha\in R$, and then we have $u_\alpha = 1$ if $\alpha\in R$ is
short, and $u_\alpha = m$ if $\alpha$ is long, where $m=2$ if $R$ is
of type $B_n$, $C_n$ or $F_4$, and $m=3$ if $R$ is of type $G_2$.

\smallskip
(iii) $R$ is not reduced, hence is of type $BC_n$
($n\geq1$). Let
$$R_1=\{\alpha\in R :\tfrac {1} {2}\alpha\notin R\} \ , \  R_2
=\{\alpha\in R : 2\alpha \notin R\},\tag 1.2$$
so that $R_1$ and $R_2$ are reduced root systems of types $B_n$, $C_n$
respectively if $n\geq 2$ (if $n=1$ they are both of type $A_1$). Up 
to similarity, there are two possibilities for $S$ when $n\geq 2$,
namely $S=R_1$ and $S=\frac {1} {2} R_2$ (which coincide when $n=1$).
In both cases $u_\alpha = 1$ or $2$ for each $\alpha\in R$ (it is for
this reason that we chose $\frac {1} {2}R_2$ rather than
$R_2$).
\medskip

Thus the function $\alpha\mapsto u_\alpha$ on $R$, when appropriately
normalized, is either constant and equal to $1$, or else takes just two
values $\{1,2\}$ or $\{1,3\}$. We shall assume this normalization
henceforth.

\remark{Remark}The classification of irreducible admissible pairs
$(R,S)$ up to similarity is closely related to (but not identical
with) the classification of irreducible affine root systems as defined
in \cite{9}, or equivalently of ``echelonnages" as defined in
\cite{3}.
\endremark 

The polynomials which are the subject of this paper will involve
parameters $q$ and $t_\alpha$, $\alpha\in R$, such that
$t_\alpha=t_\beta$ if $\vert\alpha\vert = \vert\beta\vert$. It would
be possible to regard these parameters as independent indeterminates
over $\Z$, but it will be more useful to think of them as real
variables. So let $q$ be a real number such that $0\leq q < 1$, and
for each $\alpha\in R$ let
$$q_\alpha = q^{u_\alpha},$$
so that $q_{w\alpha}=q_\alpha$ for each $w\in W$, and the set
$\{q_\alpha : \alpha\in R\}$ is either $\{q\}$ or $\{q,q^2\}$ or
$\{q,q^3\}$. From (1.1) we have 
$$q_{2\alpha}=q^2_\alpha$$
if $\alpha$, $2\alpha\in R$.

Next, for each $\alpha\in R$ let $t_\alpha$ be a real number $\geq 0$,
such that $t_\alpha=t_\beta$ if $\vert\alpha\vert=\vert\beta\vert$. If
$\alpha\in V$ but $\alpha\notin R$ we set $t_\alpha = 1$. Furthermore,
let $k_\alpha=(\log   t_\alpha)/(\log   q_\alpha)$ if $q\neq
0$ and $t_\alpha\neq 0$, so that 
$$t_\alpha=q^{k_\alpha}_\alpha.$$
If $\alpha\notin R$ we have $k_\alpha = 0$.

Finally, let $\Z[t]$ (respectively $\Z[q,t]$) denote the ring of
polynomials in the $t_\alpha$ and $t^{1/2}_{2\alpha}$ (respectively
and $q$) with integer coefficients, and let $\Q(q,t)$ denote the field
of fractions of $\Z[q,t]$, i.e., the field of rational functions of $q$
and the $t_\alpha$, $t^{1/2}_{2\alpha}$.

\head \bf \S2 \endhead

Let $(R,S)$ be an irreducible admissible pair of root systems in
$V$. Let $\{\alpha_1,\dots,\mathbreak\alpha_n\}$ be a basis (or set of simple
roots) of $R$, and let $R^+$ denote the set of positive roots
determined by this basis. Let
$$Q=\suml^n_{i=1}\Z\alpha_i, \quad \quad \quad
Q^+=\suml^n_{i=1}\N\alpha_i$$
be respectively the root lattice of $R$ and its positive octant.
Furthermore let $P$ and $P^{++}$ be respectively the weight lattice of
$R$ and the cone of dominant weights. We have $Q\subset P$ (but
$Q^+\not\subset P^{++}$ if $n > 1$). If $R$ is not reduced, then $Q$ is
the root lattice of $R_1$ (defined in (1.2)) and $P$ is the weight lattice of
$R_2$.

We define a partial order on $P$ by
$$\lambda\geq\mu \ \text{if and only} \ \lambda-\mu\in
Q^+.\tag 2.1$$

Let $A$ denote the group algebra over $\R$ of the free Abelian group
$P$. For each $\lambda\in P$, let $e^\lambda$ denote the corresponding
element of $A$, so that $e^\lambda \cdot e ^\mu = e^{\lambda + \mu}$,
$(e^\lambda)^{-1}=e^{-\lambda}$ and $e^0=1$, the identity element of
$A$. The $e^\lambda$, $\lambda\in P$, form an $\R$-basis of
$A$.

The Weyl group $W$ of $R$ acts on $P$ and hence also on
$A$: $w(e^\lambda)=e^{w\lambda}$ for $w\in W$ and $\lambda\in P$. Let
$A^W$ denote the subalgebra of $W$-invariant elements of
$A$.

Since each $W$-orbit in $P$ meets $P^{++}$ in exactly one point, it
follows that the ``monomial symmetric functions"
$$m_\lambda = \suml_{\mu\in W\lambda}e^\mu \quad \quad \quad \quad \quad (\lambda \in
P^{++})$$
form an $\R$-basis of $A^W$. Another basis is provided by the Weyl
characters: let
$$\align
R^+_2&=\{\alpha\in R^+ : 2\alpha\notin R\},\\
\rh& = \frac {1} {2}\suml_{\alpha\in R^+_2}\alpha,\tag2.2\\
\delta& =\prod_{{\alpha\in R^+_2} }
(e^{\alpha/2}-e^{-\alpha/2})=e^\rh \prod_{{\alpha\in
R^+_2}} (1-e^{-\alpha}).\tag2.3
\endalign$$
Then $w\delta = \varepsilon(w)\delta$ for each $w\in W$, where
$\varepsilon(w)=\det(w)=\pm 1$.

For each $\lambda\in P$ let
$$\chi_\lambda = \delta^{-1}\suml _{w\in W}\varepsilon(w)e^{w(\lambda +
\rh)}.\tag 2.4$$
Then $\chi_\lambda\in A^W$ for all $\lambda\in P$, and the
$\chi_\lambda$ with $\lambda\in P^{++}$ form an $\R$-basis of
$A^W$. Moreover, we have
$$\chi_\lambda = m_\lambda + \text {lower terms}\quad \quad \quad \quad (\lambda\in
P^{++})$$
where by ``lower terms" we mean a linear combination of the $m_\mu$
such that $\mu\in P^{++}$ and $\mu < \lambda$ (for the partial
ordering (2.1)).

If $\lambda\notin P^{++}$, then either $\chi_\lambda = 0$ or else there
exists $\mu\in P^{++}$ and $w\in W$ such that $\mu + \rh =
w(\lambda + \rh)$, and in this case $\chi_\lambda =
\varepsilon(w)\chi_\mu$.

Let $f\in A$, say
$$f=\suml_{\lambda\in P}f_\lambda e ^\lambda$$
with only finitely many nonzero coefficients $f_\lambda$. We shall
regard $f$ as a function on $V$ as follows: if $x\in V$, then
$$f(x)=\sum f_\lambda q ^{\langle\lambda,x\rangle}.\tag 2.5$$
For $f$ as above, define
$$\bar f= \sum f_\lambda e^{-\lambda}$$
so that $\bar f (x) = f(-x)$ for all $x\in V$. Also let
$$[f]_1=  \ \text {constant term of} \  f = f_0.$$
Clearly we have
$$[f]_1=[\bar f]_1=[wf]_1,\tag 2.6$$
for all $w\in W$.

Next, for each $\mu\in V$ we define $T_\mu f$ by
$$(T_\mu f)(x)=f(x+\mu)\tag 2.7$$
so that
$$T_\mu f = \sum f _\lambda \,q ^{\langle\lambda,\mu\rangle}
e^\lambda.$$
Each $T_\mu$ is an $\R$-algebra automorphism of $A$, with inverse
$T_{-\mu}$. We have
$$wT_\mu w^{-1}=T_{w\mu}\tag 2.8$$
for each $w\in W$, and 
$$\overline {T_\mu f} = T_{-\mu}\bar f. \tag 2.9$$

Finally, let $g=\sum g_\lambda e^\lambda$ be another element of $A$.
Then 
$$[\bar f T_\mu g]_1 = [\bar g T_\mu f]_1 .\tag 2.10$$
For both sides are equal to $\sum f_\lambda g_\lambda
q^{\langle\lambda,\mu\rangle}$.

\head \bf \S3 \endhead

We shall now define a scalar product on the algebra $A$. For this
purpose we introduce the notation
$$(x;q)_\infty=\prod _{i=0} ^{\infty}(1-xq^i)$$
and for each $k\in\R$
$$(x;q)_k=(x;q)_\infty/(x q^k ; q)_\infty.\tag 3.1$$
In particular, if $k\in\N$ we have
$$(x;q)_k=\prod _{i=0} ^{k-1}(1-x q^i). \tag 3.2$$

Now let $(R,S)$ be an irreducible admissible pair of root systems with
common Weyl group $W$, and let
$$\align\triangle =\triangle(q,t)&=\prod_{{\alpha\in R}}\frac {(t^{1/2}_{2\alpha}
e^\alpha;q_\alpha)_\infty} {(t_\alpha t^{1/2}_{2\alpha}e ^\alpha;
q_\alpha)_\infty} \tag 3.3 \\
&=\prod_{{\alpha\in R}}(t^{1/2}_{2\alpha}
e^\alpha;q_\alpha)_{k_\alpha}\endalign$$
by (3.1). If the $k_\alpha$ are all integers $\geq 0$, then by (3.2)
the product $\triangle$ is a finite product of factors of the form
$1-q^i_\alpha t^{1/2}_{2\alpha} e^\alpha$, and is clearly
$W$-invariant, hence is an element of $A^W$. In this case we define
the scalar product of two elements $f,g\in A$ to be
$$\langle f,g\rangle = \vert W\vert^{-1}[f\bar g\triangle]_1,\tag 3.4$$
i.e., the constant term of the Laurent polynomial $f\bar g\triangle$,
divided by the order of $W$.

For arbitrary values of the parameters $k_\alpha$ we proceed as
follows.
Let $Q^\lor $ be the root lattice of the dual root system $R^\lor $, and let
$T=V/Q^\lor $. Then each $e^\lambda$, $\lambda\in P$, may be regarded as a character of the torus $T$ by the rule
$e^\lambda(\overset \circ\to x)=e^{2\pi i\langle\lambda,x\rangle}$, where
$\overset \circ\to x\in T$ is the image of $x\in V$. By linearity, this
enables us to regard each element of $A$ as a continuous function on
$T$.

Consider now the product
$$(t_\alpha t^{1/2}_{2\alpha}
e^\alpha(\overset\circ\to x);q_\alpha)_\infty=\prod _{r=0}
^{\infty}\left(1-q_\alpha^{k_
\alpha + k_{2\alpha}+ r} e^\alpha(\overset\circ\to x)\right),$$
where $\alpha\in R$ and $\overset\circ\to x\in T$. This product converges
uniformly on $T$ to a continuous function (since $0\leq q_\alpha <
1$) which does not vanish on $T$ provided that $k_\alpha +
k_{2\alpha}\notin - \N$. Likewise the product
$(t^{1/2}_{2\alpha}e^\alpha;q_\alpha)_\infty$ represents a continuous
function on $T$, and therefore $\triangle$ defined by (3.3) is a
continuous function on $T$ provided that
$$k_\alpha + k _{2\alpha}\notin - \N$$
for all $\alpha\in R$ (where $k_{2\alpha}=0 \ \text {if} \
2\alpha\notin R$). Hence $\triangle$ may be expanded as a convergent
Fourier series on the torus $T$, say
$$\triangle =\sum_{\lambda\in P}a_\lambda e^\lambda, \tag 3.5$$
where
$$a_\lambda=\int_T e ^{-\lambda}\triangle,$$
the integration being with respect to normalized Haar measure on
$T$.

We now define the scalar product of $f,g\in A$ to be
$$\langle f,g\rangle = \langle f,g\rangle_{q,t} = \vert
W\vert^{-1}\int_T f\bar g\triangle.\tag3.6$$
When the $k_\alpha$ are non-negative integers, this definition agrees
with the previous one (3.4), since $\int_T e ^\lambda =
\delta_{0\lambda}$ for $\lambda\in L$.

Let
$$\triangle^+=\prod_{{\alpha\in R^+}}(t^{1/2}_{2\alpha}
e^\alpha;q_\alpha)_{k_\alpha},\tag 3.7$$
so that $\triangle = \triangle^+ \cdot \overline{\triangle^+}$. From (3.6)
it follows that
$$\langle f,g\rangle = \vert
W\vert^{-1}\int_T(f\triangle^+)(\overline{g\triangle^+})$$
and hence that the scalar product is symmetric and positive definite.

We shall next derive another expression for the scalar product (3.6)
restricted to $A^W$. For each $w\in W$ let
$$\gather
R(w)=R^+\cap - w R^+,\quad \quad \quad  t_w=\prod_{{\alpha\in
R(w)}}t_\alpha,
\tag3.8\\
W(t)=\sum_{w\in W} t_w.\tag 3.9
\endgather$$
Also let
$$\Pi =\prod_{{\alpha\in R^+}}\frac {1-t_\alpha t^{1/2}_{2\alpha} e^{-\alpha}}
{1-t^{1/2}_{2\alpha} e^{-\alpha}}$$
and
$$\triangle' = \triangle\Pi =\prod_{{\alpha\in
R^+}}(t^{1/2}_{2\alpha}e^\alpha; q_\alpha)_{k_\alpha}\,
(t^{1/2}_{2\alpha}q_\alpha e^{-\alpha}; q_\alpha)_{k_\alpha}.\tag 3.10
$$
From \cite{8} we have the identity
$$\sum_{w\in W}w\Pi = W(t)$$
so that
$$W(t)\triangle = \sum_{w\in W}w\triangle'.\tag 3.11$$

Now let $f,g\in A^W$. Then we have 
$$\align W(t)\langle f,g\rangle& = \frac {W(t)} {\vert W\vert}\int_T f
\bar g\triangle\\
&=\vert W\vert^{-1}\sum_{w\in W}\int_Tf\bar g \cdot w\triangle'\\
&=\int_Tf\bar g\triangle'
\endalign$$
since $f$ and $g$ are $W$-invariant. Hence
$$\langle f,g\rangle = W(t)^{-1}\int_T f\bar g\triangle
'\tag 3.12$$
for $f,g\in A^W$.

\remark{Remark}
If we choose to regard the parameters $q$ and $t_\alpha$ as
indeterminates over $\Z$ rather than as real numbers, we can expand
$\triangle'$ as a formal Laurent series. For this purpose let
$\ph=\sum m_i\alpha_i$ be the highest root of $R$ and let
$$x_0=q e^{-\varphi}, \ \ x_i = e^{\alpha _i}\quad (1\leq i\leq n).$$
Then $q=x_0e^\varphi = x_0x_1^{m_1}\cdots x_n^{m_n}$ is a monomial in the
$x$'s, and it follows that each of the products $q_\alpha^i e^\alpha$,
$q_\alpha^{i+1}e ^{-\alpha}$, where $\alpha\in R^+$ and $i\geq 0$, is
also a monomial in the $x$'s, since $\varphi\geq\alpha$ for each root
$\alpha\in R^+$. Moreover the total degrees of these monomials tend to
$\infty$ as $i\rightarrow\infty$, and therefore $\triangle'$ can be
expanded as a formal power series in $x_0,x_1,\dots,x_n$, say
$$\triangle' = \sum_r b_r(t)x^r,$$
where the sum is over all $r=(r_0,\dots,r_n)\in\N^{n+1}$, and
$x^r=x_0^{r_0}\cdots x_n^{r_n}$, and the coefficients $b_r(t)$ lie in
the ring $\Z[t]$ of polynomials in the $t_\alpha$ and
$t_{2\alpha}^{1/2}$ with integer coefficients.

In terms of the original variables we have
$$x^r=q^{r_0}  \exp  \(\sum^n_{i=1}r_i\alpha_i - r_0\varphi\),$$
and therefore for each $\lambda\in Q$ the coefficient of $e^\lambda$
in $\triangle'$ is
$$a'_\lambda(q,t)=\sum_{r_0}q^{r_0}b _r(t),\tag 3.13$$
where the vector $r=(r_0,r_1,\dots,r_n)$ is determined from $\lambda$
and $r_0$ by the equation
$$\sum^n_{i=1} r_i\alpha_i =\lambda + r_0\varphi,\tag 3.14$$
and the sum in (3.13) is over all integers $r_0\geq 0$ such that the
$r_i$ determined by (3.14) are all $\geq 0$. Thus we have
$$\triangle'=\sum_{\lambda\in Q} a'_\lambda(q,t)e^\lambda$$
a formal Laurent series with coefficients in
$\Z[t][[q]]$.

The identity (3.10) now gives
$$\align \triangle&=W(t)^{-1}\sum_{w\in W} w\triangle'\\
&=W(t)^{-1}\sum_{w,\lambda}a'_\lambda(q,t)e^{w\lambda}
\endalign$$
so that
$$\triangle=\sum_{\lambda\in Q}a_\lambda(q,t)e^\lambda,\tag 3.15$$
where
$$a_\lambda(q,t)=W(t)^{-1}\sum_{w\in W}a'_{w\lambda}(q,t).$$
The expression (3.15) is the expansion of $\triangle$ as a
$W$-invariant formal Laurent series, with coefficients in the ring of
formal power series $\Q(t)[[q]]$, where $\Q(t)$ is the field of
fractions of the ring $\Z[t]$. This expansion is of course the same
thing as the Fourier series (3.5).

If $f,g\in A$, the constant term in $f\bar g\triangle$ is now
well-defined, being a finite linear combination of the coefficients
$a_\lambda(q,t)$, and we have $\langle f,g\rangle = \vert
W\vert^{-1}[f\bar g\triangle]_1$ as in (3.4).
\endremark

\head \bf \S4 \endhead

We can now state the main result of this paper.

\proclaim{Theorem (4.1)} For each irreducible admissible pair $(R,S)$
of root systems there exists a unique basis $(P_\lambda)_{\lambda\in
P^{++}}$ of $A^W$ such that

\roster
\item "(i)" $P_\lambda = m_\lambda +\underset{\mu\in
P^{++}}\to{\suml_{\mu < \lambda}}
u_{\lambda\mu}(q,t)m_\mu$\newline
with coefficients $a_{\lambda\mu}(q,t)\in \Q(q,t);$
\item "(ii)" $\langle P_\lambda, P_\mu\rangle = 0  \ \text {if} \
\lambda\ne\mu$.
\endroster
\endproclaim

It is clear that the $P_\lambda$, if they exist, are unique. If the
partial order (2.1) on $P^{++}$ were a total order, the existence of
the $P_\lambda$ would follow directly from the Gram-Schmidt
orthogonalization process. However, the partial order (2.1) is not a
total order (unless rank $R=1$), and we should therefore have to
choose a compatible total order on $P^{++}$ before applying
Gram-Schmidt. The content of (4.1) is that however we extend the
partial order to a total order, we end up with the same basis of
$A^W$.

\medskip
Theorem~(4.1) will be a consequence of the following proposition:
\proclaim{Proposition (4.2) } For each irreducible admissible pair
$(R,S)$ there exists a linear operator $E:A^W\rightarrow A^W$ with the
following three properties:

\roster
\item "(i)"  E is self adjoint, i.e., $\langle Ef,g\rangle =
\langle f,Eg\rangle$  for all $f,g\in A^W$.
\item "(ii)" We have
$$Em_\lambda =\underset{\mu\in P^{++}}\to{\suml_{\mu\leq\lambda}}c_{\lambda\mu}m_\mu$$
 for each $\lambda\in P^{++}$,  with coefficients
$c_{\lambda\mu}\in q^{a(\lambda)}\Z[q,t]$, where $a:P\rightarrow 
\Q$  is a homomorphism such that $a(Q)\subset \Z$.
\item"(iii)" If $\lambda\ne\mu$, then
$c_{\lambda\lambda}\ne c_{\mu\mu}$; i.e., the eigenvalues
of $E$ are distinct.
\endroster
\endproclaim

Granted the existence of $E$ with these properties, let
$$E_\lambda = \underset{\mu\in P^{++}}\to {\prod_{{\mu <
\lambda}}}\frac {E-c_{\mu\mu}}
{c_{\lambda\lambda}-c_{\mu\mu}}$$
for $\lambda\in P^{++}$. Then the elements $P_\lambda = E_\lambda
m_\lambda$ of $A^W$ satisfy the conditions of (4.1). Indeed, it is
clear from (4.2)(ii) that the $P_\lambda$ satisfy (4.1)(i). (The
fractional exponents $a(\lambda)$ cause no trouble, because
$a(\lambda)-a(\mu)\in \Z$ if $\lambda > \mu$.) On the other hand, let
$M_\lambda$ be the subspace of $A^W$ spanned by the $m_\mu$ such that
$\mu\leq\lambda$; then $M_\lambda$ is finite-dimensional and stable 
under $E$, and the minimal polynomial of $E$ restricted to $M_\lambda$
is $\prod_{{\mu\leq\lambda}}(X-c_{\mu\mu})$, since the
$c_{\mu\mu}$ are all distinct. Hence
$(E-c_{\lambda\lambda})E_\lambda=0$ on $M_\lambda$, and therefore
$$EP_\lambda=EE_\lambda m_\lambda = c_{\lambda\lambda} E_\lambda
m_\lambda = c_{\lambda\lambda} P_\lambda.$$
If now $\lambda\ne\mu$ we have 
$$\align c_{\lambda\lambda}\langle P_\lambda,P_\mu\rangle&=\langle
EP_\lambda, P_\mu\rangle =\langle P_\lambda, EP_\mu\rangle\\
&=c_{\mu\mu}\langle P_\lambda,P_\mu\rangle
 \endalign$$
by the self-adjointness of $E$, and hence $\langle
P_\lambda,P_\mu\rangle = 0$ by (4.2)(iii).\hskip2.5cm
\vrule width3pt height10pt depth0pt

\medskip
In the next two sections we shall construct for each irreducible
admissible pair $(R,S)$ an operator $E$ satisfying the conditions of
(4.2). Our first construction, in \S5, works whenever the root system
$S^\lor $ has a minuscule fundamental weight (equivalent conditions are
that $P\ne Q$, or that $R$ is not of type $E_8$, $F_4$ or $G_2$).
In \S6 we shall give another construction which workes in
these excluded cases.

\head \bf \S5 \endhead

In this section we shall assume that $S^\lor $ possesses a minuscule
fundamental weight, i.e., that there exists a vector $\pi\in  V$ such
that $\langle \pi,\alpha_\ast\rangle$ takes just two values $0$ and
$1$ as $\alpha$ runs through $R^+$. We have then (2.7)
$$T_\pi e^\alpha = q^{\langle\pi,\alpha\rangle} e^\alpha =
q_\alpha^{\langle\pi,\alpha_\ast\rangle} e^\alpha$$
so that 
$$
\aligned 
T_\pi e^\alpha =\cases  q_\alpha e^\alpha  &\text{if }\langle\pi,\alpha_\ast\rangle = 1,\\
e^\alpha &  \text{if } \langle\pi,\alpha_\ast\rangle = 0.
\endcases\endaligned\tag 5.1$$
Now let 
$$\align
\Phi_\pi &=(T_\pi\triangle^+)/\triangle^+\\
&=\underset{\langle\pi,\alpha_\ast\rangle=1}\to{\prod_{{\alpha\in
R^+}}}\frac {1-t_\alpha t^{1/2}_{2\alpha}
e^\alpha} {1-t^{1/2}_{2\alpha}e^\alpha}\tag5.2
\endalign$$ 
by (5.1) and the definition (3.7) of $\triangle^+$.

We define an operator $E_\pi$  on $A$ as follows:
$$E_\pi f = \sum_{w\in W} w(\Phi_\pi\cdot T_\pi f).\tag5.3$$
Let us first show that $E_\pi$ is self-adjoint (on the assumption that
it maps $A$ into $A$, which we shall justify shortly). Since
$$E_\pi f = \sum_{w\in W}\frac {w(T_\pi(\triangle^+ f))} {w\triangle^+},$$
and since $\triangle = w\triangle =
w\triangle^+\cdot \overline{w\triangle^+}$ for each $w\in W$, we have
$$\align 
\langle E_\pi f, g\rangle &=\vert W\vert^{-1}\sum_{w\in W}
\left[w(T_\pi(\triangle^+ f))\cdot\overline{w(\triangle^+ g)}\right]_1\\
&=\left[T_\pi(\triangle^+f)\cdot \overline{\triangle^+ g}\right]_1
\endalign$$
by (2.6), and by (2.10) this expression is symmetrical in $f$ and $g$.
Hence
$$\langle E_\pi f, g\rangle = \langle E_\pi g, f\rangle = \langle f,
E_\pi g\rangle.\tag5.4$$

To show that $E_\pi$ maps $A$ into $A$, we need to express $\Phi_\pi$
in a more convenient form. If $\alpha\in R^+$ and $\frac {1}
{2}\alpha\in R^+$, the corresponding factors in the product (5.2)
combine to give
$$\align
\frac {(1-t_\alpha e^\alpha)\left(1-t_{\alpha/2}t^{1/2}_\alpha e
^{\alpha/2}\right)}
{(1-e^\alpha)\left(1-t_\alpha^{1/2}e^{\alpha/2}\right)}&=\frac
{\left(1+t_\alpha^{1/2}e^{\alpha/2}\right)\left(1-t_{\alpha/2}t^{1/2}_\alpha
e^{\alpha/2}\right)} {1-e^\alpha}\\
&=\frac {1+(1-t_{\alpha/2})t_\alpha^{1/2}
e^{\alpha/2}-t_{\alpha/2}t_\alpha e^\alpha} {1-e^\alpha}.
\endalign$$
If $\frac {1} {2}\alpha\notin R^+$ (and $2\alpha\notin R^+$), this is
still correct, since then $t_{\alpha/2}=1$. It follows that
$$\Phi_\pi=\prod_{{\alpha\in R^+_2}}\frac
{1+\left(1-t_{\alpha/2}^{\langle\pi,\alpha_\ast\rangle}\right)
t_\alpha^{\langle\pi
,\alpha_\ast\rangle/2} e^{\alpha/2}-
(t_{\alpha/2}t_\alpha)^{\langle\pi,\alpha_\ast\rangle}e^\alpha}
{1-e^\alpha},\tag 5.5$$
where as before $R^+_2=\{\alpha\in R^+:2\alpha\notin R\}$.

Since
$t_\alpha^{\langle\pi,\alpha_\ast\rangle} =
q^{k_\alpha\langle\pi,\alpha\rangle}$, we have
$$\prod_{{\alpha\in R^+}}
t_\alpha^{\langle\pi,\alpha_\ast\rangle} = q ^{2\langle\pi,\rh_
k\rangle},\tag5.6$$
where
$$\rh_k=\frac {1} {2}\sum_{\alpha\in R^+} k_\alpha\alpha. \tag5.7$$
Hence the product (5.5) for $\Phi_\pi$ can be rewritten in the form
$$\Phi_\pi=\delta^{-1}e^\rho q^{2\langle\pi,\rh _k\rangle}\Psi,$$
where $\delta$ and $\rh$ are as in (2.2) and (2.3), and $\Psi$ is the
product
$$\Psi=\prod_{{\alpha\in
R^+_2}}\left(1+
\left(t_{\alpha/2}^{-\langle\pi,\alpha_\ast\rangle}-1\right)t_
\alpha^{-\langle\pi,\alpha_\ast\rangle/2}e^{-\alpha/2}-(t_{\alpha/2}
t_\alpha)^{-\langle\pi,\alpha\ast\rangle}e ^{-\alpha}\right).$$
If we multiply out this product, we shall obtain
$$\Phi_\pi=\delta^{-1}q^{2\langle\pi,\rh_
k\rangle}\sum_X \varphi_X(t)e^{\rh-\sigma(X)},\tag5.8$$
summed over all subsets $X$ of $R^+$ such that
$$\alpha\in X\Rightarrow 2\alpha\notin X, \tag 5.9$$
with the following notation:
$$\align
\sigma(x)&=\sum_{\alpha\in X}\alpha,\\
\varphi_X(t)&=\prod_{{\alpha\in
X}}\varphi_\alpha(t),\tag 5.10
\endalign$$
$$\aligned\varphi_\alpha(t)=\cases
-(t_{\alpha/2}t_\alpha)^{-\langle\pi,\alpha_\ast\rangle}&\text
{if} \ 2\alpha\notin R,\\
\left(t_\alpha^{-\langle\pi,\alpha_\ast\rangle}-1\right) 
\, t_{2\alpha}^{-\langle\pi,
\alpha_\ast\rangle/2}&\text {if} \ 2\alpha\in
R.\endcases\endaligned\tag 5.11$$
We can now calculate $E_\pi e^\mu$, where $\mu\in P$. Since $w(T_\pi e
^\mu)=q^{\langle\pi,\mu\rangle} e^{w\mu}$, and since
$w\delta=\varepsilon(w)\delta$ for each $w\in W$, where $\varepsilon(w)= 
\det(w)=\pm 1$, we shall obtain from (5.8)
$$\align 
E_\pi e^\mu &= \delta^{-1} q^{\langle\pi, 2\rh _k +
\mu\rangle}\sum_X\varphi_X(t)\sum_{w\in
W}\varepsilon(w)e^{w(\mu+\rh-\sigma(X))}\\
&=q^{\langle\pi,2 \rh_ k +\mu\rangle}\sum_X\varphi_X(t)\chi_{\mu
-\sigma(X)}
\endalign$$
from which it follows that $E_\pi$ maps $A$ into $A^W$.

Now let $\lambda\in P^{++}$. Then we have
$$\align E_\pi m_\lambda&= \sum_{\mu\in W\lambda} E_\pi e^\mu\\
&=\sum_X\varphi_X(t)\sum_{\mu\in W\lambda} q ^{\langle\pi,2 \rho_k +
\mu\rangle}\chi_{\mu-\sigma(X)}.\tag 5.12
\endalign$$
In this sum, either $\chi_{\mu-\sigma(X)}= 0$ or else there exists
$w\in W$ and $\nu\in P^{++}$ such that
$$\nu + \rh = w(\mu + \rh -\sigma(X))\tag 5.13$$
in which case $\chi_{\mu - \sigma(X)}=\varepsilon(w)\chi_\nu$. But $\rh -
\sigma(X)$ is of the form
$$\rh - \sigma(X)=\frac {1} {2}\sum_{\alpha\in R_2^+} 
\varepsilon_\alpha\alpha,$$
where each coefficient $\varepsilon_\alpha$ is $\pm 1$ or $0$, hence
$w(\rh - \sigma(X))$ is of the same form, and therefore
$$w(\rh - \sigma(X))=\rh - \sigma(Y)\tag 5.14$$
for same subset $Y$ of $R^+$ such that $\alpha\in Y\Rightarrow
2\alpha\notin Y$. From (5.13) and (5.14) it follows that
$$\nu=w\mu -\sigma(Y)\leq w\mu\leq\lambda\tag 5.15$$
and hence that $E_\pi m_\lambda$ is a linear combination of the
$\chi_\nu$ such that $\nu\in P^{++}$ and $\nu\leq\lambda$. Hence we
have
$$E_\pi m_\lambda = \underset{\nu\in P^{++}}\to{\sum_{\nu\leq\lambda}}
b_{\lambda\nu}\chi_\nu,$$
where from (5.12) the coefficient $b_{\lambda\nu}$ is given by
$$b_{\lambda\nu} = \sum\varepsilon(w)q^{\langle\pi, 2\rh_k
+\mu\rangle}\varphi_X(t)$$
summed over triples $(X,\mu,w)$ where $X \subset R^+$, $\mu\in W\lambda$ and
$w\in W$ satisfy (5.9) and (5.13). From (5.6) and the definition
(5.10) of $\varphi _X(t)$ it follows that
$q^{\langle\pi,2\rh_k\rangle}\varphi_X(t)\in \Z[t]$, the ring of
polynomials over $\Z$ generated by the $t_\alpha$ and
$t_{2\alpha}^{1/2}$. As to the scalar product $\langle\pi,\mu\rangle$,
we have
$$\langle\pi,\mu\rangle=\langle\pi,w_0\lambda\rangle +
\langle\pi,\theta\rangle$$
where $w_0$ is the longest element of $W$ and $\theta
=\mu-w_0\lambda\in Q^+$. Now $\langle\pi,\alpha_\ast\rangle = 0$ or 1
for each $\alpha\in R^+$, and hence $\langle\pi,\alpha\rangle = 0$ or
$u_\alpha$ for $\alpha\in R^+$. Hence $\langle\pi,\theta\rangle$ is a
non-negative integer and therefore $b_{\lambda\nu}\in
q^{\langle\pi,w_0\lambda\rangle}\Z[q,t]$. The exponent
$\langle\pi,w_0\lambda\rangle$ need not be an integer, but it is a
rational number.

From this it follows that
$$E_\pi m_\lambda = \underset{\nu\in P^{++}}\to{\sum_{\nu\leq\lambda}}
c_{\lambda\nu}(\pi)m_\nu\tag 5.16$$
with coefficients $c_{\lambda\nu}(\pi)\in
q^{\langle\pi,w_0\lambda\rangle}\Z[q,t]$.

We must now calculate the leading coefficient
$c_{\lambda\lambda}(\pi)$ in (5.16). From (5.15) it follows that
$\nu=\lambda$ if and only if $Y$ is empty and $w\mu =\lambda$, that is
to say if and only if $\mu=w^{-1}\lambda$ and $w(\rh
-\sigma(X))=\rh$, or equivalently $\sigma(X)=\rh - w^{-1}\rh$. But
this implies \cite{8} that $X=R_2(w)=R_2^+\cap-wR_2^+$. Hence the
coefficient of $m_\lambda$ in (5.16) is
$$c_{\lambda\lambda}(\pi)=\sum_{w\in W}\varepsilon(w)\varphi_ {R_2
(w)}q^{\langle\pi,2\rh_k  + w^{-1}\lambda\rangle}.$$
Since $\varepsilon(w)=(-1)^{\vert R_2(w)\vert}$, we obtain from (5.10) and
(5.11)
$$\align\varepsilon(w)\varphi_{R_2(w)}(t)&=\prod_{{\alpha\in
R_2(w)}}(t_{\alpha/2} t_\alpha)^{-\langle\pi,\alpha_\ast\rangle}\\
&=\prod_{\alpha\in R(w)}
t_\alpha^{-\langle\pi,\alpha_\ast\rangle}\\
&=q^{\langle\pi,w^{-1}\rh_k -\rh_k\rangle}
\endalign$$
since
$t_\alpha^{-\langle\pi,\alpha_\ast\rangle}=q^{-\langle\pi,k_\alpha\alpha\rangle}$. Hence
$$\align 
c_{\lambda\lambda}(\pi)&=q^{\langle\pi,\rh _k\rangle}\sum_{w\in W}q
^{\langle w\pi,\lambda + \rh_k\rangle}\tag 5.17\\
&=q^{\langle\pi, \rh _k\rangle} \widetilde m_\pi(\lambda +\rh _k)
\endalign$$
where $\widetilde m_\pi =\suml_{w\in W} e^{w\pi}$.

It remains to examine whether the eigenvalues
$c_{\lambda\lambda}(\pi)$ of $E_\pi$ are all distinct as $\lambda$
runs through $P^{++}$, for a suitable choice of the minuscule weight
$\pi$. It will appear that this is so in all cases except $D_n$, $n\geq
4$ (and of course excepting $E_8,F_4$ and $G_2$, where there is no
minuscule weight).

Let
$$p_r(x)=\sum_{w\in W}\langle x,w\pi\rangle^r\quad \quad \quad (x\in
V).$$
The $p_r$, $r\geq 1$, are $W$-invariant polynomial functions on $V$.

\proclaim{(5.18)} Suppose that $S$ is not of type $D_n$ {\rm(}$n\geq
4${\rm)}, and
that if $S$ is of type $A_n$ the minuscule weight $\pi$ is the
fundamental weight corresponding to an end  node of the Dynkin
diagram. Then the $p_r$ generate the $\R$-algebra of $W$-invariant
polynomial functions on $V$, and hence separate the $W$-orbits in
$V$.
\endproclaim

This is easily verified for $S$ of type $A,B$ or $C$. For $E_6$ and
$E_7$ see \cite{12}.

\medskip
Assume now that the hypotheses of (5.18) are satisfied, and that
$\lambda,\mu\in P^{++}$ are such that
$c_{\lambda\lambda}(\pi)=c_{\mu\mu}(\pi)$, i.e., that
$$\sum_{w\in W} q^{\langle\lambda +\rh_k,w\pi\rangle} = \sum_{w\in W}
q^{\langle\mu + \rh_k, w\pi\rangle}.$$
By operating on both sides with $(q\partial/\partial q)^r$ and then
setting $q=1$ and $t_\alpha=1$ for each $\alpha\in R$, 
we obtain $p_r(\lambda)=p_r(\mu)$ for all
$r\geq1$. Hence by (5.18) $\lambda$ and $\mu$ are in the same
$W$-orbit, and therefore $\lambda=\mu$. It follows that the
eigenvalues $c_{\lambda\lambda}(\pi)$ of $E_\pi$ are all
distinct.

There remains the case where $S(=R)$ is of type $D_n$. Let
$\varepsilon,\dots,\varepsilon_n$ be an orthonormal basis of $V$; we may
then assume that $R^+$ consists of the vectors $\varepsilon_i\pm\varepsilon
_j$ with $i < j$. Then $P^{++}$ consists of the vectors
$\lambda=\sum\lambda_i\varepsilon_i$ for which the $\lambda_i$ are all
integers or all half-integers, and
$\lambda_1\geq\cdots\geq\lambda_{n-1}\geq\vert\lambda_n\vert$. The
fundamental weights
$$\pi_1 = \frac {1} {2}(\varepsilon_1 +\cdots + \varepsilon_n), 
\quad \pi_2 = \pi_1
-\varepsilon_n$$
are both minuscule, and the $W$-orbit of $\pi_1$ (respectively $\pi_2$)
consists of all sums $\frac {1} {2}\sum\pm\varepsilon_i$ containing an
even (respectively odd) number  of minus signs. Hence the formula (5.17)
gives
$$c_{\lambda\lambda}(\pi_1)\pm c_{\lambda\lambda}(\pi_2)= n! \prodl
_{i=1} ^{n}\left(q^{\lambda_ i}t^{n-i}\pm q ^{-\lambda _i}\right)$$
where $t=t_\alpha,\alpha\in R$.

Choose an integer $N > \frac {1} {2} n(n-1)$. Then the eigenvalues of
the operator
$$E=\frac {1} {n!}\left(\left(t^N + 1\right)E_{\pi _1}+ 
\left(t^N - 1\right)E_{\pi_2}\right)\tag 5.19$$
are from above
$$c_\lambda=t^N\prodl _{i=1} ^{n}(q^{\lambda _i} t^{n-i} + q
^{-\lambda_i}) +\prodl _{i=1}
^{n} (q^{\lambda_i} t^{n-i} - q ^{-\lambda_i}).$$
Now suppose that $\lambda,\mu\in P^{++}$ are such that $c_\lambda =
c_\mu$. From our choice of $N$ it follows that
$$\prodl _{i=1} ^{n}(q^{\lambda_i} t ^{n-i} + q^{-\lambda _i})=\prodl
_{i=1} ^{n}(q^{\mu_i} t ^{n-i} + q^{-\mu_i}), \tag 5.20$$
$$\prodl _{i=1} ^{n}(q^{\lambda_i} t ^{n-i} - q^{-\lambda_i})=\prodl
_{i=1} ^{n}(q^{\mu_i} t ^{n-i} - q ^{-\mu_i}).\tag 5.21$$
From (5.20) we conclude that $\lambda_i =\mu_i$ ($1\leq i\leq n-1$) and
$\lambda_n =\pm\mu_n$.
If $\lambda_n\ne 0$, then (5.21) shows that $\lambda_n=\mu_n$, and
hence $\lambda=\mu$. So the eigenvalues $c_\lambda$ of $E$ are all
distinct.

To recapitulate, let $E:A^W\rightarrow A^W$ be the operator $E_\pi$
defined by (5.3) when $R$ is not of type $D_n$, and by (5.19) when $R$
is of type $D_n$. By (5.4), (5.16) and the discussion above, it
follows the $E$ satisfies the three conditions of (4.2).

\head \bf \S6 \endhead

In the cases where there is no minuscule weight available, another
construction is needed. We shall assume in this section that $R$ is
reduced, so that either $S=R$ or $S=R^\lor $, from  the classification in
\S1.

Let $\varphi\in R^+$ be such that $\varphi_\ast$ is the highest root
of $S$, and let $\pi=(\varphi_\ast)^\lor =u_\varphi\varphi ^\lor $. For each
$\alpha\in R^+$ the Cauchy-Schwarz inequality gives
$$0\leq\langle\pi,\alpha_\ast\rangle = \frac
{2\langle\varphi_\ast,\alpha_\ast\rangle} {\vert\varphi_\ast\vert^2}\leq
\frac {2\vert\alpha_\ast\vert} {\vert\varphi_\ast\vert}\leq 2$$
with equality if and only if $\alpha_\ast = \varphi_\ast$. Since
$\langle\pi,\alpha_\ast\rangle$ is an integer, it follows that for
each $\alpha\in R^+$
$$\langle\pi,\alpha_\ast\rangle = \cases 0 \ \text {or} \  1 \ &\text  
{if} \
\alpha\ne\varphi,\\
2 \ &\text   {if} \ \alpha = \varphi.\endcases\tag6.1$$
Thus $\pi$ just fails to be a minuscule weight.

\remark{Remark}
In fact $\varphi_\ast =\varphi$, so that $u_\varphi =1$ and
$q_\varphi=q$. This is clear if $S=R$, whereas if $S=R^\lor $ (and $R$ has
two root-lengths) $\varphi$ is the highest short root of $R$, so that
$u_\varphi = 1$ in this case also. Hence $\pi=\varphi^\lor $.
\endremark

Let
$$\Phi_\pi = (T_\pi\triangle^+)/\triangle^+$$
as in \S5, and define an operator $F_\pi$ on $A^W$ as follows:
$$F_\pi f =\suml_{w\in W} w(\Phi_\pi \cdot U_\pi f)$$
where $U_\pi = T_\pi -1$.

Let us first show that $F_\pi$ is self-adjoint. If $f,g\in A^W$ we
have 
$$\align\langle F_\pi f,g\rangle &= \vert W\vert^{-1}\suml_{w\in
W}\left[w((T_\pi\triangle^+)(U_\pi
f))\cdot\overline{w(\triangle^+g)}\right]_1\\
&=\left[(T_\pi\triangle^+)(T_\pi f - f)\overline{\triangle^+g}\right]_1\\
&=\left[T_\pi(\triangle^+ f)\cdot\overline{\triangle^+ g}\right]_1 -
\left[(T_\pi\triangle^+)\overline{\triangle^+} \cdot f \bar g\right]_1\\
&=A(f,g)-B(f,g),
\endalign$$
say. We have $A(f,g)=A(g,f)$ by (2.10). As to $B(f,g)$, let
$G=(T_\pi\triangle^+)\overline{\triangle^+}$ and let $w_0$ be the longest
element of the Weyl group $W$. Then $w_0\pi = -\pi$, and
$w_0\triangle^+ =\overline{\triangle^+}$, so that
$$w_0 G = \big(T_{-\pi}\overline{\triangle^+}\big)\triangle^+
=\overline{T_\pi\triangle^+}\triangle^+ = \overline G$$
and therefore
$$\align B(f,g)&=[G f\bar g]_1 =[(w_0 G)f\bar g]_1\\
&=[\overline G f\bar g]_1=[G\bar f g]_1 = B(g,f).
\endalign$$
It follows that
$$\langle F_\pi f,g\rangle = \langle F_\pi g,f\rangle = \langle f,
F_\pi g\rangle\tag 6.2$$
and hence that $F_\pi$ is self-adjoint.

From (6.1) we have, since $q_\varphi = q$,
$$\align \Phi_\pi &= \frac{1-t_\varphi e^\varphi}{1-e^\varphi}
\cdot\frac {1-q t_\varphi e^\varphi} {1-q
e^\varphi}\underset{\alpha\ne\varphi}\to{\prod_{{\alpha\in
R^+}}}\frac {1-t_\alpha^{\langle\pi,\alpha_\ast\rangle}e^\alpha}
{1-e^\alpha}\\
&=q^{2\langle\pi,\rh _k\rangle}\frac
{(1-t_\varphi^{-1}e^{-\varphi})(1-q^{-1}t_\varphi^{-1}e^{-\varphi})}
{(1-e^{-\varphi})(1-q^{-1}e^{-\varphi})}\underset{\alpha\ne\varphi}
\to{\prod_{{\alpha\in R^+}}}\frac
{1-t_\alpha^{-\langle\pi,\alpha_\ast\rangle}e^{-\alpha}} {1-e^{-\alpha}}\\
&=\delta^{-1}e^\rh q^{2\langle\pi,\rh _k\rangle}\frac {(1-t^{-1}_\varphi
e^{-\varphi})(1-q^{-1}t_ \varphi^{-1}e^{-\varphi})}
{1-q^{-1}e^{-\varphi}}\underset{\alpha\ne\varphi}\to{\prod_{
{\alpha\in R^+}}}
\left(1-t_\alpha^{-\langle\pi,\alpha_\ast\rangle}e^{-\alpha}\right)
\endalign$$
and therefore
$$\Phi_\pi=\delta^{-1}q^{2\langle\pi,\rh _k\rangle}\frac {(1-t^{-1}_\varphi
e^{-\varphi})(1-q^{-1}t_\varphi^{-1} e^{-\varphi})}
{1-q^{-1}e^{-\varphi}}\suml_X\psi_X(t)e^{\rh-\sigma(X)}\tag6.3$$
summed over all subsets $X$ of $R^+$ such that $\varphi\notin X$, where
$$\sigma(X)=\suml_{\alpha\in X}\alpha,$$
and
$$\psi_X(t)=(-1)^{\vert X\vert}\prod_{{\alpha\in X}}t_\alpha
^{-\langle\pi,\alpha_\ast\rangle}.\tag6.4$$

Now let $\lambda\in P^{++}$ and $\mu\in W\lambda$. Since $T_\pi e^\mu
= q^{\langle\mu,\pi\rangle}e^\mu=q^{\langle\mu,\varphi^\lor \rangle} e^\mu$, we
have
$$U_\pi m_\lambda = \suml_{\mu\in
W\lambda}\left(q^{\langle\mu,\varphi^\lor \rangle}-1\right)e^\mu.\tag6.5$$
Any $\mu\in W\lambda$ such that $\langle\mu,\varphi^\lor \rangle = 0$
will contribute nothing to this sum. The remaining elements of the orbit
$W\lambda$ fall into pairs $\{\mu,w_\varphi\mu\}$ where
$\langle\mu,\varphi^\lor \rangle > 0$ and $w_\varphi$ is the reflection
associated with $\varphi$. Hence we may rewrite (6.5) in the form
$$U_\pi m_\lambda = \underset{\langle\mu,\varphi^\lor \rangle >
0}\to{\underset{\mu\in W\lambda}\to\suml}
\left(q^{\langle\mu,\varphi^\lor \rangle}-1\right)e^\mu\left(1-(q
e^\varphi)^{-\langle\mu,\varphi^\lor \rangle}\right)$$
from which it follows that
$$\frac {U_\pi m_\lambda}
{1-q^{-1}e^{-\varphi}}=\underset{\langle\mu,\varphi^\lor \rangle >
0}\to{\underset{\mu\in
W\lambda}\to\suml}
\left(q^{\langle\mu,\varphi^\lor \rangle}-1\right)\suml _{j=0}
^{\langle\mu,\varphi^\lor \rangle -1}q^{-j}e^{\mu -j\varphi}.\tag 6.6$$
From (6.3) and (6.6) we obtain
$$\multline
\Phi_\pi\cdot U_\pi m_\lambda = \delta^{-1}q^{2\langle\pi,\rh_
k\rangle}\suml_{X,\mu}\psi_X(t)e^{\rh-\sigma(X)}(1-t_\varphi^{-1}
e^{-\varphi})\\
\times
(1- q^{-1}t_\varphi^{-1}
e^{-\varphi})
\left(q^{\langle\mu,\varphi^\lor \rangle}-1\right)\prodl _{j=0} ^{\langle
\mu,\varphi^\lor \rangle-1}q^{-j}e^{\mu - j\varphi}.
\endmultline\tag 6.7$$
This is a sum of terms of the form $a\delta^{-1}e^\eta$, $\eta \in P$.
Since $\suml_{w\in W}w(\delta^{-1}e^\eta)\in A^W$, it follows from
(6.7) that $F_\pi m_\lambda\in A^W$, and hence that $F_\pi$ maps $A^W$
into $A^W$.

Moreover, the terms $a\delta^{-1}e^\eta$ that occur in (6.7) are such
that
$$\eta=\rh-\sigma(X)+\mu-j\varphi$$
where $0\leq j\leq \langle\mu,\varphi^\lor \rangle +1$ and
$\langle\mu,\varphi^\lor \rangle \geq1$ and $X\subset R^+-\{\varphi\}$.
If $\eta$ is not regular (i.e., if $W_\eta\ne 1$, where $W_\eta$ is
the subgroup of $W$ that fixes $\eta$) it will contribute
nothing to (6.7). If on the other hand $\eta$ is regular, then we have
$w\eta = \xi + \rh$ for some $\xi \in P^{++}$ and some
$w\in W$, so that
$$\xi + \rh = w(\rh -\sigma(X))+w(\mu-j\varphi).\tag6.8$$
There are two cases to consider.

\medskip\noindent
(i) Suppose that $w\varphi=\alpha\in R^+$. Since $\rh-\sigma(X)$ is of
the form $\frac {1} {2}\suml_{\alpha\in R^+}\varepsilon_\alpha\alpha$,
where each $\varepsilon_\alpha$ is $\pm 1$, it follows that
$w(\rh-\sigma(X))$ is of the same form, hence that
$$w(\rh-\sigma(X))=\rh-\sigma(Y)$$
for some subset $Y$ of $R^+$. Hence
$$\xi=w\mu-j\alpha-\sigma(Y)\leq w\mu\leq\lambda\tag6.9$$
and therefore each such term $a\delta^{-1}e^\eta$ in (6.7) contributes
$\varepsilon(w)a\chi_\xi$, where $\xi\leq\lambda$, to $F_\pi
m_\lambda$.

Moreover we have equality in (6.9) if and only if $j=0$, the subset
$Y$ is empty, and $w\mu=\lambda$, i.e., $\mu=w^{-1}\lambda$ and
$\sigma(X)=\rh-w^{-1}\rh$, so that $X=R(w)=R^+\cap-wR^+$. The
coefficient $a$ of $\delta^{-1}e^\eta=\delta^{-1}e^{w^{-1}(\lambda +
\rh)}$ in (6.7)  is then
$$a=\psi_{R(w)}(t)q^{2\langle\pi,\rh_
k\rangle}\left(q^{\langle\lambda,w\pi\rangle}-1\right)\tag6.10$$
(since $q_\varphi^{\langle\mu,\varphi^\lor \rangle} = q^{\langle
w^{-1}\lambda,\pi\rangle}=q^{\langle\lambda, w\pi\rangle}$).
From (6.4) we have
$$\psi_{R(w)}(t)=\varepsilon(w)\prod_{{\alpha\in
R(w)}}t_\alpha^{-\langle\pi,\alpha_\ast\rangle}$$
and since
$t_\alpha^{-\langle\pi,\alpha_\ast\rangle}=q^{-\langle\pi,k_\alpha
\alpha\rangle}$ it follows that
$$\psi_{R(w)}(t)=\varepsilon(w)q^{-\langle\pi,\rh _k - w^{-1}\rh_
k\rangle}.\tag 6.11$$
From (6.10) and (6.11) the coefficient of $\delta^{-1}e^{w^{-1}(\lambda +
\rh)}$ in (6.7) is therefore
$$a=\varepsilon(w)q^{\langle\pi,\rh _k\rangle}\left(q^{\langle w\pi,\lambda
+ \rh _k\rangle}-q^{\langle w\pi,\rh _k\rangle}\right)\tag 6.12$$
and the corresponding contribution to $F_\pi m_\lambda$ is
$$q^{\langle\pi,\rh _k\rangle}\left(q^{\langle w\pi,\lambda +\rh
_k\rangle}-q^{\langle w\pi,\rh_k\rangle}\right)\chi_\lambda.\tag6.13$$

\smallskip\noindent
(ii) Suppose now that $w\varphi = -\alpha$, where $\alpha\in R^+$, and
let $\nu=w_\varphi\mu=\mu-\langle\mu,\varphi^\lor \rangle\varphi$. Then
$$\mu-j\varphi=\nu-\varphi +j'\varphi,$$
where $j'=\langle\mu,\varphi^\lor \rangle + 1 - j$, so that $0\leq
j'\leq\langle\mu,\varphi^\lor \rangle + 1$. Hence (6.8) now takes the form
$$\xi + \rh = w(\rh - \sigma(X)-\varphi)+ w(\nu + j'\varphi).$$
Since $\varphi\notin X$, we have
$w(\rh-\sigma(X)-\varphi)=\rh-\sigma(Y)$ for some $Y\subset R^+$,
and therefore
$$\xi = w\nu -\sigma(Y)-j'\alpha\leq w\nu\leq\lambda.\tag6.14$$
So again each term $a\delta^{-1}e^\eta$ in (6.7) contributes
$\varepsilon(w)a\chi_\xi$, where $\xi\leq\lambda$, to $F_\pi
m_\lambda$.

Moreover, we have equality in (6.14) if and only if $j'=0$, the subset
$Y$ is empty, and $w\nu=\lambda$, i.e., $\nu=w^{-1}\lambda$ and
$\sigma(X)+\varphi =\rh-w^{-1}\rh$, so that $X\cup\{\varphi\}=R(w)$ and
$\mu=w_\varphi w^{-1}\lambda$ and $j=\langle\mu,\varphi^\lor \rangle + 1$.
Hence the coefficient of $\delta^{-1}e^\eta =\delta^{-1}e^{w^{ -1}(\lambda
+ \rh)}$ in $\Phi_\pi\cdot U_\pi m_\lambda$ is now
$$a=\psi_{R(w)-\{\varphi\}}(t)q^{2\langle\pi,\rh
_k\rangle}t_\varphi^{-2}\cdot q^{-\langle\mu,\varphi^\lor
\rangle}\left(q^{\langle
\mu,\varphi^\lor \rangle}-1\right).$$
Since $-\langle\mu,\varphi^\lor \rangle = -\langle w_\varphi
w^{-1}\lambda,\varphi^\lor \rangle =\langle w^{-1}\lambda,\varphi^\lor \rangle
= \langle\lambda,w\varphi^\lor \rangle$
we have
$$a=-t^{-2}_\varphi\psi
_{R(w)-\{\varphi\}}(t)q^{2\langle\pi,\rh_k\rangle}\left(q^{\langle\lambda,w
\pi\rangle}-1\right).$$
Moreover, from (6.4),
$$\align
-t_\varphi^{-2}\psi_{R(w)-\{\varphi\}}(t)&=\varepsilon(w)t^{-2}_\varphi
\underset{\alpha\ne\varphi}\to{\prod_{{\alpha\in R(w)}}}
t_\alpha^{-\langle\pi,\alpha_\ast\rangle}\\
&=\varepsilon(w)\prod_{{\alpha\in R(w)}}
t_\alpha^{-\langle\pi,\alpha_\ast\rangle}\\
&=\varepsilon(w)q^{-\langle\pi,\rh_k-w^{-1}\rh_k\rangle}
\endalign$$
as in (6.11). So finally the coefficient of
$\delta^{-1}e^{w^{-1}(\lambda +\rh)}$ in (6.7) is given by the same
expression (6.12) as before, and the corresponding contribution to
$F_\pi m_\lambda$ is again given by (6.13).

To recapitulate, these calculations show that $F_\pi m_\lambda$ is a
linear combination of the Weyl characters ${\chi_\xi}$ such that $\xi\in
P^{++}$ and $\xi\leq\lambda$, and an inspection of (6.7) shows that
the coefficient of $\chi_\xi$ in $F_\pi m_\lambda$ lies in
$q^{-\langle\lambda,\varphi^\lor \rangle}\Z[q,t]$. Moreover the
coefficient of $\chi_\lambda$ is
$$q^{\langle\pi,\rh _k\rangle}\suml_{w\in W}\left(q^{\langle w\pi,\lambda
+\rh _k\rangle}-q^{\langle w\pi,\rh _k\rangle}\right)
=q^{\langle\pi,\rh _k\rangle}(\widetilde m _\pi(\lambda +\rh_k)-\widetilde
m_\pi(\rh_k))
\tag 6.15$$
where $\widetilde m_\pi=\suml_{w\in W}e ^{w\pi}=\vert W_\pi\vert m_\pi$.
Hence we have
$$F_\pi m_\lambda =\underset{\mu\in
P^{++}}\to{\underset{\mu\leq\lambda}\to\suml}c_{\lambda\mu}(\pi)m_\mu
\tag 6.16$$
with $c_{\lambda\lambda}(\pi)$ given by (6.15), and
$c_{\lambda\mu}(\pi)\in q^{-\langle\lambda,\varphi^\lor \rangle}\Z[q,t]$.

To complete the proof of (4.2), it remains to establish that the
$c_{\lambda\lambda}(\pi),\lambda\in P^{++}$, are all distinct when $R$
is of type $E_8$, $F_4$, or $G_2$. For this purpose we argue as in the
last part of \S5: if $\lambda,\mu\in P^{++}$ are such that
$c_{\lambda\lambda}(\pi)=c_{\mu\mu}(\pi)$, then by operating with
$(q\partial/\partial q)^r$ and then setting $q=1$ we shall obtain
$$\suml_{w\in W}\langle w\pi,\lambda\rangle^r=\suml_{w\in W}\langle
w\pi,\mu\rangle^r$$
for each $r\geq 1$. But in each case $\pi$ is either the highest
root or the highest short root of $R$, and it is known \cite{12} that
when $R$ is of type $E_8$, $F_4$ or $G_2$ the polynomial functions on
$V$
$$p_r(x)=\suml_{w\in W}\langle w\pi, x\rangle^r\quad \quad \quad (r\geq
1)$$
generate the $\R$-algebra of $W$-invariant polynomial functions in
$V$, and therefore separate the $W$-orbits in $V$. Hence $\lambda$ and
$\mu$ lie in the same $W$-orbit, and so $\lambda=\mu$.

By (6.2), (6.16) and the above discussion, it follows that the linear
operator $F_\pi$ on $A^W$ satisfies the three conditions of (4.2) when
$R$ is of type $E_8$, $F_4$ or $G_2$. This completes the proof of
(4.2) and hence of Theorem~(4.1).

\head \bf \S7 \endhead

Let us say that a linear operator $L:A^W\rightarrow A^W$ is
{\it triangular} if
$$Lm_\lambda = \alpha_\lambda m_\lambda +\text {lower terms}$$
for each $\lambda\in P^{++}$, that is to say if the matrix of $L$
relative to the basis $(m_\lambda)$ of $A^W$ is triangular.
\proclaim{(7.1)} If $L$ is triangular and self-adjoint, the $P_\lambda$ are
eigenfunctions of $L$.
\endproclaim

\demo{Proof} Since $L$ is triangular, we have
$$LP_\lambda = \underset{\mu\in
P^{++}}\to{\underset{\mu\leq\lambda}\to\suml}\alpha_{\lambda\mu}P_\mu$$
with coefficients $a_{\lambda\mu}$ given by
$$a_{\lambda\mu}\vert P_\mu\vert^2 = \langle L P_\lambda, P_\mu\rangle
= \langle P_\lambda, L P_\mu\rangle$$
since $L$ is self-adjoint. But
$$LP_\mu =\underset{\nu\in P^{++}}\to{\underset{\nu\leq\mu}\to\suml}
a_{\mu\nu} P_\nu,$$
so that $\langle P_\lambda, LP_\mu\rangle = 0$ unless $\mu =
\lambda$. Hence $a_{\lambda\mu}=0$ unless $\mu =\lambda$, which proves
(7.1).
\enddemo

From (7.1) it follows that all self-adjoint triangular linear
operators on $A^W$ are simultaneously diagonalized by the $P_\lambda$,
and therefore commute with each other.

Consider in particular the case in which $R$ is of type $A_n$. Then
all the fundamental weights $\pi_i$ ($1\leq i\leq n$) are minuscule, and
hence the construction of \S5 furnishes $n$ self-adjoint triangular
linear operators $E_{\pi_ i}$ ($1\leq i\leq n$) on $A^W$, which by the
above remark commute with each other. Moreover the eigenvalues of
$E_{\pi_ i}$ are (up to a scalar factor) $m_{\pi_ i}(\lambda + \rh _k)$.
Since the $m_{\pi_ i}$ generate the algebra $A^W$, it follows that for
each $f\in A^W$ there is a unique linear operator $L_f\in \R[E_{\pi
_1},\dots,E_{\pi_n}]$ such that
$$L_fP_\lambda = f(\lambda + \rh_k)P_\lambda\tag 7.2$$
for all $\lambda\in P^{++}$. Moreover, since each $E_{\pi_ i}$ is a
linear combination of the translation operators $T_{w\pi_i}$, $w\in
W$, it follows that $L_f$ is of the form
$$L_f=\suml_{\mu\in P}a_{\mu,f}(q,t)T_\mu\tag7.3$$
with coefficients $a_{\mu,f}(q,t)$ which are rational functions of $q$
and $t$.

For an arbitrary admissible pair $(R,S)$ of root systems, we can of
course define $L_f$ for each $f\in A^W$ by (7.2). But except in the
case $A_n$ just mentioned, I do not know whether $L_f$ is expressible
in the form (7.3).

\head \bf \S8 \endhead

In the following sections we shall consider some particular
cases.

\medskip\noindent
(i) Suppose first that $k_\alpha = 0$ for each $\alpha\in R$, i.e.,
$t_\alpha = 1$. Then $\triangle = 1$, so that the scalar product is now
$$\langle f,g\rangle =\vert W\vert^{-1}[f\bar g]_1.$$
It follows that
$$P_\lambda = m_\lambda\tag 8.1$$
for all $\lambda\in P^{++}$, and hence that
$$\langle P_\lambda,P_\mu\rangle = \vert
W_\lambda\vert^{-1}\delta_{\lambda\mu}\tag 8.2$$
where $W_\lambda$ is the subgroup of $W$ that fixes $\lambda$.

\smallskip\noindent
(ii) Next, suppose that $R$ is reduced and that $k_\alpha = 1$ for
each $\alpha\in R$, so that $t_\alpha = q_\alpha$. Then
$$\triangle=\prod_{{\alpha\in R}}
(1-e^\alpha)=\delta\bar\delta$$
and the scalar product is now
$$\langle f, g\rangle = \vert W\vert^{-1}[f\delta \cdot \overline {g \delta}]_1.$$
Hence in this case we have
$$P_\lambda = \chi_\lambda\tag8.3$$
for all $\lambda\in P^{++}$, and
$$\langle P_\lambda, P_\mu\rangle =\delta_{\lambda\mu}.\tag 8.4$$

\head \bf \S9 \endhead

In this section we shall consider the cases where $R$ has rank $1$,
hence is of type $A_1$ or $BC_1$.

\medskip\noindent
(i) Suppose first that $R$ is of type $BC_1$. In this case it terms out
that the polynomials $P_\lambda$ are a particular case of the
orthogonal polynomials defined by Askey and Wilson \cite{2}. They
define
$$\multline
p_n(x;a,b,c,d\mid q)\\
=a^{-n}(ab;q)_n\,(ac;q)_n\,(ad;q)_n \cdot 
{}_4\varphi_3\!\left[\matrix  q^{-n},q^{n-1}abcd,ae^{i\theta},
ae^{-i\theta}\\ ab,ac,ad  \endmatrix; q,q\right]
\endmultline\tag9.1$$
where $x= \cos \theta$, the parameters $a,b,c,d$ and $q$
lie in the interval $(-1, 1)$ of $\R$, and $_4\varphi _3$ is the usual
notation for a $q$-hypergeometric series. Thus $p_n$ is a Laurent
polynomial in $e^{i\theta}$ with leading term $(abcd q^{n-1}; q)_n
e^{ni\theta}$, and in fact is (despite its appearance) symmetrical in
all four parameters $a,b,c,d$. For our purposes it is more convenient
to take $e^{i\theta}$ rather than $\cos\theta$ as the argument, and
we prefer to have the leading coefficient equal to $1$; so we define
$$\tilde p_n(e^{i\theta};a,b,c,d\mid q)=(abcd q^{n-1};q)^{-1}_n\,
p_n(\cos\theta;a,b,c,d\mid q).\tag9.2$$
Then Theorem~(2.2) of \cite{2} takes the form
$$\frac {1} {2\pi}\int^\pi_0\tilde p_ m \,\tilde p_ n\left\vert
f(e^{i\theta})\right\vert^2 d\theta = \delta_{mn}\tilde h_n\tag 9.3$$
where
$$f(u)=(u^2;q)_\infty/(au;q)_\infty(bu;q)_\infty(cu;q)_\infty(du;q)
_\infty\tag 9.4$$
and
$$\tilde h_n = \frac {(A q^{2n-1};q)_\infty\,(A q^{2n};q)_\infty}
{(A q^{n-1};q)_\infty \,(q^{n+1};q)_\infty}\,\Pi\,(a b
q^n;q)^{-1}_\infty\tag 9.5$$
where $A=a b c d$ and $\Pi$ means the product of the six terms such
as $(ab q^n; q)^{-1}_\infty$.

Now let $R$ be a root system of type $BC_1$, with $R^+=\{\alpha,
2\alpha\}$ and $S=\{\pm\alpha\}$. Then $q_\alpha = q$ and
$q_{2\alpha}=q^2$. We shall write $k_1,k_2, t_1,t_2$ for
$k_\alpha,k_{2\alpha}, t_\alpha, t_{2\alpha}$ respectively, so
that $t_1=q^{k_1}$ and $t_2=q^{2k_2}$. We have
$P^{++}=\N\alpha$.

Let
$$(a, b, c, d)=(q^{1/2}, - q^{1/2}, t_1t_2^{1/2}, - t_2^{1/2}).\tag 9.6$$
With this choice of parameters, $f(u)$ takes the form
$$f(u)=\frac {(t^{1/2}_2 u; q)_\infty\,(u^2;q^2)_\infty} {(t_1t_2^{1/2}u;
q)_\infty\,(t_2 u^2;q^2)_\infty}$$
so that $\triangle = f(e^\alpha)f (e^{-\alpha})$. It now follows from
(9.1), (9.2) and (9.3) that for $n\geq 0$
$$P_{n\alpha}=\tilde p_n(e^\alpha;q^{1/2},-q^{1/2},t_1t_2^{1/2},-t
_2^{1/2}\mid q)\tag 9.7$$
and that $\vert P_{n\alpha}\vert^2$ is given by (9.5) when the
parameters $a,b,c,d$ are as in (9.6). After some reduction we find that
$$\vert P_{n\alpha}\vert^2=\frac
{(q^{2n}t_1t_2;q)_{k_1}\,(q^{2n}t^2_1t_2;q^2)_{k_2}}
{(q^{2n+1}t_2;q)_{k_1}\,(q^{2n+2};q^2)_{k_2}}\tag 9.8$$

For later reference we shall calculate $P_{n\alpha}
(\rh^\ast_k)$, where $\rh^\ast_k=\frac {1} {2}(k_1 + k_2)\alpha^\lor $,
so that $e^\alpha(\rh^\ast_k)=q^{k_1 + k_2} = t_1t_2^{1/2}$. The
original formula (9.1) for $p_n$ shows that when $e^{i\theta}$ is
replaced by $a$, the series $_4\varphi_3$ reduces to its first term,
which is 1. By the symmetry of the parameters $a, b, c, d$ it
follows that when $e^{i\theta}$ is replaced by $c$ in $p_n$, the
result is
$$c^{-n}(ac;q)_n\,(bc;q)_n\,(cd;q)_n.$$
In view of (9.6) and (9.7), this observation enables us to calculate
$P_{n\alpha}(\rh_k^\ast)$; after some reduction, we obtain
$$P_{n\alpha}(\rh_k^\ast)=(t_1t_2^{1/2})^{-n}\frac
{(q^{2n}t_1t_2;q)_{k_1}\,(q^{2n}t_1^2t_2;q^2)_{k_2}}
{(t_1t_2;q)_{k_1}\,(t^2_1 t_2;
q^2)_{k_2}}.\tag 9.9$$

\smallskip\noindent
(ii) Suppose now that $R$ is of type $A_1$, with positive root
$\alpha$, so that $P^{++}=\frac {1} {2}\N\alpha$. In this case the
$P_\lambda$ are essentially the $q$-ultraspherical polynomials of
Askey and Ismail \cite{1}. If $\lambda =\frac {1} {2}n\alpha$, where
$n\geq 0$, we have
$$P_\lambda=\frac {(q;q)_n} {(t; q)_n}\varphi_n(e^{\alpha/2})\tag 9.10$$
where $t=t_\alpha$ and 
$$\varphi_n(x)=\suml_{i + j = n}\frac {(t;q)_i\,(t;q)_j}
{(q;q)_i\,(q;q)_j}x^{i-j}\tag9.11$$
(so that $\varphi_n(e^{i\theta})=C_n(\cos\theta; t\mid q)$ in the
notation of \cite{1}). The $\varphi_n$ have the generating function
$$F(x,u)=\suml_{n\geq 0}\varphi_n(x)u^n=1/(xu;q)_k(x^{-1}u;q)_k\tag
9.12$$
(where $k=k_\alpha$, so that $t=q^k$), as follows from (9.11) and the
$q$-binomial theorem.

One way of establishing (9.10) is to verify that the $\varphi_n$
defined by (9.11) are eigenfunctions of the operator $E_\pi$ of \S5,
which in the present situation takes the form
$$(E_\pi f)(x)=\frac {tx - x^{-1}} {x-x^{-1}}f(q^{1/2}x)+\frac
{x-tx^{-1}} {x-x^{-1}}f(q^{1/2}x^{-1})\tag 9.13$$
for a Laurent polynomial $f(x)\in \R [x,x^{-1}]$. More precisely, we
have to verify that
$$E_\pi\varphi_n = (tq^{n/2}+ q^{-n/2})\varphi_n$$
or equivalently (9.12) that
$$(E_\pi F)(x,u)=t F(x,q^{1/2}u)+F(x,q^{-1/2}u).$$
But this is straightforward to verify from (9.12) and (9.13).

From \cite{1} we have
$$\vert P_{n\alpha/2}\vert^2=(q^n t; q)_k/(q^{n+1};q)_k\tag 9.14$$
and
$$P_{n\alpha/2}(\tfrac {1} {2}k\alpha^\lor )=t^{-n/2}(q^n
t;q)_k/(t;q)_k;\tag 9.15$$
this latter formula is equivalent to
$\varphi_n(t^{1/2})=t^{-n/2}(t^2;q)_n/(q;q)_n$, which follows from
(9.12) since
$$\align 
F(t^{1/2},ut^{1/2})&=1/(u;q)_k(tu;q)_k\\
&=1/(u;q)_{2k}\\
&=\suml_{n\geq 0}\frac {(t^2;q)_n} {(q;q)_n} u^n.
\endalign$$

\head \bf \S10 \endhead

We consider next the case where $q=0$, the $t_\alpha$ being
arbitrary. Then we have
$$\triangle =\prod_{{\alpha\in R}}\frac
{1-t_{2\alpha}^{1/2}e^\alpha} {1-t_{2\alpha}^{1/2}t_\alpha
e^\alpha}.$$
In this case there is an explicit formula for $P_\lambda$:
$$P_\lambda = W_\lambda(t)^{-1}\suml_{w\in W}
w\left(e^\lambda\prod_{{\alpha\in R^+}}\frac {1-t_\alpha
t_{2\alpha}^{1/2}e^{-\alpha}} {1-t_{2\alpha}^{1/2}e^{-\alpha}}\right)\tag
10.1$$
{\it where $W_\lambda$ is the subgroup of $W$ that fixes $\lambda\in
P^{++}$, and}
$$W_\lambda(t)=\suml_{w\in W_\lambda} t_w$$
{\it with $t_w$ as defined in} (3.8).

To prove (10.1), we shall first show that $P_\lambda = m_\lambda$ +
lower terms, and then that $\langle P_\lambda, P_\mu\rangle = 0  \ \text
{if} \ \lambda\ne\mu$.
 
Let
$$\align
\Phi_\lambda &= 
e^\lambda\prod_{{\alpha\in R^+}}\frac
{1-t_\alpha t_{2\alpha}^{1/2}e^{-\alpha}}
{1-t_{2\alpha}^{1/2}e^{-\alpha}}\\
&=\delta^{-1}e^{\lambda + \rh}\prod_{{\alpha\in
R^+_2}}\left(1+(1-t_{\alpha/2})t_\alpha^{1/2}e^{-\alpha/2}-t_{\alpha
/2}t_\alpha e^{-\alpha}\right)
\endalign$$
as in \S5, where $R^+_2=\{\alpha\in R^+ : 2\alpha\notin R\}$ and
$\delta$, $\rh$ are as defined in (2.2), (2.3). On multiplying out this
product we shall obtain
$$\Phi_\lambda = \delta^{-1}\sum_X\varphi_X(t)e^{\lambda + \rh
-\sigma(X)}\tag 10.2$$
summed over all subsets $X$ of $R^+$ such that $\alpha\in X\Rightarrow
2\alpha\notin X$, where $\sigma(X)= \suml_{\alpha\in X}\alpha$ and
$$\align\varphi_X(t)&=\prod_{{\alpha\in
X}}\varphi_\alpha(t),\tag 10.3\\
\varphi_\alpha(t)&=\cases \quad \  -t_{\alpha/2} t_\alpha \ &\text {if} \
2\alpha\notin R,\\
(1-t_\alpha)t_{2\alpha}^{1/2} \ &\text {if} \ 2\alpha\in R.\endcases\endalign$$
Let
$$Q_\lambda = \sum_{w\in W} w\Phi_\lambda.$$
Then it follows from (10.2) that
$$Q_\lambda = \sum_X\varphi_X(t)\chi_{\lambda-\sigma(X)}\tag10.4$$
summed over subsets $X\subset R^+$ as above. If
$\chi_{\lambda-\sigma(X)}\ne 0$, there exists $w\in W$ and $\mu\in
P^{++}$ such that
$$\mu + \rh = w(\lambda + \rh - \sigma(X))$$
and we have $\chi_{\lambda-\sigma(X)}=\varepsilon(w)\chi_\mu$. Now (5.14)
$w(\rh - \sigma(X))=\rh - \sigma(Y)$ for some subset $Y$ of $R^+$ such
that $\alpha\in Y\Rightarrow 2\alpha\notin Y$. Hence $\mu = w\lambda -
\sigma(Y)\leq w\lambda\leq\lambda$, and so it follows from (10.4) that
$Q_\lambda$ is a linear combination of the $\chi_\mu$ such that $\mu\in
P^{++}$ and $\mu\leq\lambda$.

Moreover, we have $\mu=\lambda$ if and only if $w\lambda=\lambda$ and
$\sigma(Y)=0$, that is to say if and only if $w\in W_\lambda$ and
$\sigma(X)=\rh -w^{-1}\rh$, which implies that $X=R_2(w)$ and hence
(10.3)
$$\align \varphi_X(t)&=\prod_{{\alpha\in
R_2(w)}}(-t_{\alpha/2}t_\alpha)\\
&=\varepsilon(w)\prod_{{\alpha\in R(w)}}t_\alpha =\varepsilon(w)t_w.
\endalign$$
Hence the coefficient of $\chi_\lambda$ in $Q_\lambda$ is
$$\suml_{w\in W_\lambda} t_w = W_\lambda(t)$$
and therefore $P_\lambda$ as defined by (10.1) is of the form
$m_\lambda$ + lower terms.

It remains to prove that $\langle P_\lambda, P_\mu\rangle = 0 \ \text
{if} \ \lambda\ne\mu$. We may assume that $\lambda\nleq\mu$. We
have
$$Q_\lambda\triangle = \sum_{w\in W} w\Psi_\lambda$$
where
$$\align \Psi_\lambda&=\triangle\Phi_\lambda = e^\lambda\triangle ^+ =
e^\lambda\prod_{{\alpha\in R^+}}\frac
{1-t_{2\alpha}^{1/2}e^\alpha} {1-t_\alpha t_{2\alpha}^{1/2}e^\alpha}\\
&=e^\lambda\prod_{{\alpha\in R^+}}\bigg(1+(t_\alpha -1)\sum_{r\geq
1}t^{r-1}_\alpha t_{2\alpha}^{r/2}e^{r\alpha}\bigg)\\
&=\sum_{\mu\in Q^+}a_\mu e^{\lambda + \mu}
\endalign$$
say, with $a_0 = 1$. Hence
$$Q_\lambda\triangle = \sum_{\mu\in Q^+}a_\mu\sum_{w\in W}e^{w(\lambda
+ \mu)}.$$
If $\lambda + \mu = w_1\pi$, where $\pi\in P^{++}$ and $w_1\in W$, then
we have $\pi\geq w_1\pi =\lambda + \mu\geq\lambda$, with equality only
if $\mu=0$. Hence
$$Q_\lambda\triangle =\vert W_\lambda\vert m_\lambda + \text {higher
terms};$$
and since $\lambda\nleq\mu$ this sum has no terms in common with
$$P_\mu = m_\mu + \text {lower terms}.$$
Hence $\langle Q_\lambda, P_\mu\rangle =\vert W\vert^{-1}[\bar {P}_\mu
Q_\lambda\triangle]_1 =0$. This completes the proof of (10.1).

Moreover, these calulations show that
$$\langle Q_\lambda,P_\lambda\rangle = \frac {\vert W_\lambda\vert}
{\vert W\vert}[m_\lambda\overline {m_\lambda}]_1 = 1$$
so that
$$\vert P_\lambda\vert^2 = W_\lambda (t)^{-1}.\tag 10.5$$

The formula (10.1) is essentially the formula of \cite{7}, 
Theorem~(4.1.2) for the zonal spherical function on a $p$-adic Lie group. More
precisely, let $G$ be a simply-connected group of $p$-adic type, as
defined in \cite{7}, and let $K$ be a special maximal compact
subgroup of $G$, such that the root system $\Sigma_1$ of \cite{7}, (3.1)
is the {\it dual} $R^\lor $ of $R$. The root structure of $G$ attaches a
positive integer $q_{\alpha^\lor}$ to each root $\alpha^\lor \in R^\lor$, 
and we
take  $t_\alpha = q_{\alpha^\lor}^{-1}$. The double cosets of $K$ in $G$ are
indexed by the elements of $Q^{++}=P^{++}\cap Q$, and the zonal spherical
functions $\omega_S$ on $G$ relative to $K$ are parametrized by the
$\C$-algebra homomorphisms $s:\C[Q]^W\rightarrow\C$.

For each $\lambda\in Q^{++}$ let $g_\lambda$ be a representative of
the corresponding double coset of $K$ in $G$. Then the formula for the
zonal spherical function is
$$\omega_s(g^{-1}_\lambda)=u_\lambda(t)s(P_\lambda)$$
with $P_\lambda$ as in (10.1) and
$$u_\lambda(t)=\frac {W_\lambda (t)} {W(t)}\prod_{{\alpha\in
R^+}}t_\alpha^{\langle\lambda,\alpha^\lor \rangle/2}.$$
Moreover, $\triangle$ is essentially the Plancherel measure on the
space of positive definite zonel spherical functions on $G$ relative
to $K$.

\head \bf \S11 \endhead

In this section we shall consider the ``limiting case" as
$q\rightarrow 1,$ the parameters $k_\alpha$ remaining fixed. We shall
assume that
$$k_\alpha\geq 0\tag 11.1$$
for all $\alpha\in R$.\newline

Let
$$\triangle_k =\prod_{{\alpha\in R}}(1-e^\alpha)^{k_\alpha}$$
considered (as in \S3) as a continuous function on the torus
$T=V/Q^\lor .$ For $f,g\in A$ we define
$$\langle f,g\rangle_k=\vert W\vert^{-1}\int _{T} ^{} f\bar
{g}\triangle_k\tag 11.2$$
using $\triangle_k$ in place of $\triangle(q,t)$. As before, this
scalar product on $A$ is symmetric and positive definite.

Suppose that, in addition to (11.1), we have
$$k_\alpha + 2k_{2\alpha}\geq 1\tag 11.3$$
for all $\alpha\in R$ (so that $k_\alpha\geq1$ if
$2\alpha\notin R$). Then
$$\lim_{q\rightarrow 1}\triangle(q,t)=\triangle_k\tag 11.4$$
uniformly on $T$.

This is a consequence of the following fact \cite{6}: if $r,s\in \R$
and $z\in\C$ then 
$$\lim_{{q\rightarrow 1}}\frac {(q^rz;q)_\infty}
{(q^sz;q)_\infty}=(1-z)^{s-r}$$
uniformly on the disc $\vert z\vert\leq 1$, provided that $r\leq s$
and $r+s\geq 1$. If we take $q=q_\alpha$, $r=k_{2\alpha}$, $s=k_\alpha
+ k_{2\alpha}$ we obtain
$$\lim_{{q\rightarrow
1}}(t_{2\alpha}^{1/2}e^\alpha;q_\alpha)_{k_\alpha}
=(1-e^\alpha)^{k_\alpha}$$ 
uniformly on $ T$, provided that (11.1) and
(11.3) hold. Taking the product over all $\alpha\in R$, we obtain
(11.4).

Until further notice we shall assume (11.3) as well as (11.1).\newline

Let
$$f(q)=\suml_{\lambda\in P}f_\lambda(q)e^\lambda$$
be an element of $A$ depending on $q\in (0,1)$. If
$f_\lambda(q)\rightarrow f_\lambda$ as $q\rightarrow 1$ for each
$\lambda\in P$, we shall write
$$\lim_{{q\rightarrow 1}} f(q)=f$$
where $f=\sum f_\lambda e^\lambda$.

Suppose also that
$$\lim_{{q\rightarrow 1}}g(q)=g$$
in $A$. Then
$$\lim_{{q\rightarrow 1}}\langle f(q),
g(q)\rangle_{q,t}=\langle f,g\rangle_k.\tag11.5$$

By linearity it is enough to prove this when
$f(q)=f_\lambda(q)e^\lambda$, and $g(q)=g_\mu(q)e^
\mu$. We have then
$$\langle f(q), g(q)\rangle _{q,t}=\vert
W\vert^{-1}f_\lambda(q)g_\mu(q)\int_T e^{\lambda -\mu}\triangle(q,t)$$
which by (11.4) tends to the limit
$$\vert W\vert^{-1}f_\lambda g_\mu\int_T e^{\lambda - \mu}\triangle_k
=\langle f,g\rangle_k$$
as $q\rightarrow 1$.

Consider now the behaviour of
$$P_\lambda(q,t)=\underset{\mu\in
P^{++}}\to{\underset{\mu\leq\lambda}\to\suml}u_{\lambda\mu}
(q,t)m_\mu\quad \quad \quad\quad \quad \quad  (\lambda\in P^{++})\tag11.6$$
as $q\rightarrow 1$. We claim that
$$\lim_{{q\rightarrow 1}} P_\lambda (q,t) \ \text {\it exists for
each} \  \lambda\in P^{++}:\tag11.7$$
in other words, that each of the coefficients $u_{\lambda\mu}(q,t)$
tends to a finite limit as $q\rightarrow 1$.

We shall prove this by induction on $\lambda$. When $\lambda =0$
there is nothing to prove, since $P_0(q,t)=1$; so assume that
$\lambda\ne 0$ and that
$$P_\mu(k)=\lim_{{q\rightarrow 1}}P_\mu(q,t)\tag 11.8$$
exists for all $\mu\in P^{++}$ such that $\mu < \lambda$.

The equations (11.6) can be inverted to give say
$$m_\lambda =\underset{\mu\in
P^{++}}\to{\underset{\mu\leq\lambda}\to\suml}v_{\lambda\mu}
(q,t)P_\mu(q,t)\tag 11.9$$
with $v_{\lambda\lambda}=1$. The $v$'s are cofactors of the
(unipotent) matrix formed by the $u$'s, hence are polynomials in the
$u$'s, and conversely the $u$'s are polynomials in the $v$'s. From
(11.9) and the orthogonality of the $P$'s we have
$$v_{\lambda\mu}(q,t)=\frac {\langle m_\lambda, P_\mu(q,t)\rangle
_{q,t}} {\vert P_\mu(q,t)\vert ^2_{q,t}}.$$
Hence, by (11.5) and (11.8), we have
$$\lim_{{q\rightarrow 1}}v_{\lambda\mu}(q,t)=\frac
{\langle m_\lambda, P_\mu(k)\rangle_k} {\vert P_\mu(k)\vert^2_k}$$
whenever $\mu < \lambda$. But, as we have just remarked,
$u_{\lambda\mu}(q,t)$ is a polynomial in the $v$'s with integer
coefficients. Hence
$$u_{\lambda\mu}(k)=\lim_{{q\rightarrow
1}}u_{\lambda\mu}(q,t)$$
exists for each $\mu <\lambda$, and (11.7) is proved.

Since $u_{\lambda\mu}(q,t)\in \C(q,t)$, hence is a rational function
of $q$ and the $q^{k_\alpha}$, its limit as $q\rightarrow 1$ may be
computed by differentiating its numerator and denominator sufficiently
often, and then setting $q=1$. This shows that $u_{\lambda\mu}(k)$ is
rational function of the $k_\alpha$.

We now define
$$\align 
P_\lambda(k)&=\suml_{\mu\leq\lambda}u_{\lambda\mu}(k)m_\mu\tag 11.10\\
&=\lim_{{q\rightarrow 1}}P_\lambda (q,t).
\endalign$$
By (11.5) we have
$$\langle P_\lambda(k),P_\mu(k)\rangle_k=0\tag11.11$$
if $\lambda\ne\mu$, and the properties (11.10) and (11.11)
characterize the polynomials $P_\lambda(k)$.

The existence of these polynomials has been established by Heckman
and Opdam \cite{4, 14} by other methods.

Define linear operators $\square$ and $D_\alpha(\alpha\in R)$ on $A$
by
$$\square e^\lambda =\vert\lambda\vert^2 e^\lambda , \quad D_\alpha
e^\lambda=\langle\lambda,\alpha\rangle e^\lambda$$
and as in \cite{14} let
$$L(k)=\square + \frac {1} {2}\suml_{\alpha\in R}k_\alpha\frac
{1+e^{-\alpha}} {1-e^{-\alpha}}D_\alpha.\tag11.12$$
An equivalent definition is
$$L(k)f=\delta_k^{-1/2}\square\left(\delta^{1/2}_k
f\right)-\left(\delta_k^{-1/2}\square\delta^{1/2}_k\right)f\tag11.13$$
where
$$\delta_k=\prod_{{\alpha\in
R^+}}(e^{\alpha/2}-e^{-\alpha/2})^{2k_\alpha}.$$
From (11.12), it is clear that $wL(k)w^{-1}=L(k)$ for all $w\in W$,
and that if $f\in A$ then $g=\delta L (k)f\in A$, where $\delta$ is as
defined in \S2. If $f\in A^W$, it follows that $g$ is
skew-symmetric with respect to $W$, and hence $L(k)f=\delta^{-1}g\in
A^W$. Thus $L(k)$ maps $A^W$ into $A^W$, and a simple calculation
shows that
$$L(k)m_\lambda =\langle\lambda,\lambda + 2\rh _k\rangle m_\lambda 
 \text { + lower terms}$$
for $\lambda\in P^{++}$, where $\rh_k$ is given by (5.7). Moreover,
since $\triangle = \delta^{1/2}_k\overline{\delta^{1/2}_k}$, it follows from
(11.13) that $L(k)$ is self-adjoint for the scalar product (11.2).
Hence as in (7.1) we conclude that the $P_\lambda$ are eigenfunctions
of $L(k)$, and more precisely that
$$L(k)P_\lambda(k)=\langle\lambda,\lambda + 2\rh_k\rangle P_\lambda
(k)\tag 11.14$$
for $\lambda\in P^{++}$.

From (11.12) and (11.14) we derive the following recurrence relation
for the coefficients $u_{\lambda\mu}(k)$ in (11.10):
$$  (\vert\lambda +\rh_k\vert^2-\vert\mu
+\rh_k\vert^2)u_{\lambda\mu}(k)
=2\suml_{\alpha\in R^+}\suml_{r\geq 1}k_\alpha\langle \mu +
r\alpha,\alpha\rangle u_{\lambda,\mu + r\alpha}(k).
\tag11.15$$
Here $\lambda,\mu\in P^{++}$, $\mu < \lambda$, and
$u_{\lambda\mu}=u_{\lambda,w\mu}$ for all $w\in W$.
Let $\nu =\lambda - \mu\in Q^+$, then
$$\vert\lambda + \rh_k\vert^2 - \vert\mu
+\rh_k\vert^2=\langle\nu,\lambda +\mu + 2\rh_k\rangle.$$
We have $\langle\nu,\lambda + \mu\rangle > 0$; also 
$$\langle\rh_k,\alpha^\lor _i\rangle = k_{\alpha_i}+2k_{2\alpha_i}$$
for a simple root $\alpha_i\in R$, from which it follows that
$\langle\rh_k,\nu\rangle$ is a positive linear combination of the
$k_\alpha$. Hence if the $k_\alpha$ are all $\geq 0$, as we are
assuming throughout (11.1), we have
$$\vert\lambda + \rh_k\vert^2 -\vert\mu +\rh_k\vert^2 > 0$$
whenever $\lambda,\mu\in P^{++}$ and $\lambda >\mu$. It follows now
from the recurrence formula (11.15) by induction on $\nu=\lambda -\mu$
that the coefficients $u_{\lambda\mu}(k)$ are {\it positive}, and more
precisely that they are rational functions of the $k$'s in which both
numerator and denominator are polynomials in the $k$'s with positive
integral coefficients. So finally we can drop the restriction (11.3);
the polynomials $P_\lambda$ are well-defined provided that (11.1)
holds.

Now let $G/K$ be a non-compact symmetric space, and $G=KAN$ an Iwasawa
decomposition of the semisimple Lie group $G$. Let $\frak g,\frak k,\frak
a$ be the
Lie algebras of $G,K,A$ respectively, let $\frak a^*$ be the vector space
dual to $\frak a$, and $\Sigma\subset \frak a^*$ the restricted root system of
$G/K$. For each $\beta\in\Sigma$ let $m_\beta$ denote the multiplicity
of $\beta$. We shall take
$$V=\frak a^*, \quad R=2\Sigma,\quad  
k_\alpha=\tfrac {1} {2}m_{\alpha/2}\tag 11.16$$
for each $\alpha\in R$. Then $L(k)$ is the radial part of the Laplacian
on $G/K$.

Let $\frak h\supset \frak a$ be a Cartan subalgebra of $\frak g$; 
then $\frak h=\frak a\oplus
(\frak h\cap \frak k)$, and we put $\frak h_\R =\frak  a\oplus i(\frak h
\cap\frak  k)$. Let $M$ be a finite-dimensional irreducible representation
space for $G$, with highest weight $\lambda\in \frak h^*_\R$. Then
(\cite{5}, p.~535) $M$ is spherical, i.e., has a nonzero vector fixed by
$K$, if and only if $\lambda$ vanishes on $i(\frak h\cap \frak k)$ 
and (regarded
as an element of $V=\frak a^*$) we have
$\langle\lambda,\alpha^\lor \rangle\in\N$ for all $\alpha\in R^+$, i.e.,
$\lambda\in P^{++}$.

The recursion formula (11.15) conincides with Harish-Chandra's recursion
formula (\cite{5}, p.~427) for the coefficients of the zonal spherical
function $\omega_\lambda$ defined by $M$. Hence for $X\in \frak a^+$, the
positive Weyl chamber in $\frak a$, we have
$$\omega_\lambda(\exp X)=P_\lambda(X)/P_\lambda(0)\tag11.17$$
where each $e^\lambda$ ($\lambda\in P$) is now to be regarded as a
function on $\frak a=V^*$ by the rule $e^\lambda(X)=e^{\lambda(X)}$, the
latter $e$ being the classical exponential function.

If $G/K$ is compact, the formula (11.17) is the same, except that now
$X\in i\frak a$ in place of $\frak a^+$.

Note the contrast with the $p$-adic situation of \S10. There the
restricted root system was $R^\lor $, and the weights $\lambda\in P^{++}$
indexed the double cosets of $K$ in $G$. Here, on the other hand, the
restricted root system is (similar to) $R$, and the $\lambda\in P^{++}$
index the finite-dimensional spherical representations.

\head \bf \S12 \endhead

We shall conclude with some conjectures. In order to state them
concisely, we introduce the $q$-gamma function, defined for $0 < q < 1$
by
$$\Gamma _q(x)=\frac {(q;q)_{x-1}}  {(1-q)^{x-1}}=\frac
{(q;q)_\infty\,
(1-q)^{1-x}} {(q^x;q)_\infty},$$
where $x\in \R$ but $x\ne 0,-1,-2,\dots$, at which points
$\Gamma_q$ has simple poles.\newline
We have
$$\Gamma_q(x + 1) = \frac {1-q^x} {1-q}\Gamma_q(x)$$
and hence
$$\lim_{{k\rightarrow 0}} \ \frac {\Gamma_q(kx)}
{\Gamma_q(ky)}=\lim_{{k\rightarrow 0}}  \ \frac {1-q^{ky}}
{1-q^{kx}}\frac {\Gamma_q(kx + 1)} {\Gamma_q(ky + 1)} = \frac {y}
{x}.\tag12.1$$

As $q\rightarrow 1$, $\Gamma_q(x)\rightarrow\Gamma(x)$ for each $x$,
and therefore we shall write $\Gamma_1(x)$ for the ordinary gamma
function.

It will be convenient also to introduce
$$\Gamma_q^*(x)=1/\Gamma_q(1-x).\tag 12.2$$
When $q=1$ we have
$$\Gamma_1^*(x)=\frac {\sin \pi x} {\pi}\Gamma_1(x)\tag12.3$$
but there is no particularly simple relationship between
$\Gamma_q^*(x)$ and $\Gamma_q(x)$ for general values of $q$.

Now define, for $\lambda\in V$ and $\alpha\in R^+$,
$$c_\alpha(\lambda;q_\alpha) \ = \ \frac
{\Gamma_{q_\alpha}\!\!\left(\langle\lambda,\alpha^\lor \rangle + \frac {1}
{2}k_{\alpha/2}\right)}
{\Gamma_{q_\alpha}\!\!\left(\langle\lambda,\alpha^\lor \rangle+ \frac {1}
{2}k_{\alpha/2}+ k_\alpha\right)}\tag12.4$$
(where as usual $k_{\alpha/2}=0 \ \text {if} \ \frac {1}
{2}\alpha\notin R$), and 
$$c(\lambda)=c(\lambda;q,t)=\prod_{{\alpha\in R^+}}
c_\alpha(\lambda;q_\alpha).\tag 12.5$$
Also define $c^*_\alpha(\lambda;q_\alpha)$ and $c^*(\lambda)$ by using
$\Gamma^*$ in place of $\Gamma$ in (12.4) and (12.5).

We can now state
\proclaim{Conjecture (12.6)} For all $\lambda\in P^{++}$
$$\vert P_\lambda\vert^2 = \frac {c^\ast(-\lambda -\rh_k)} {c(\lambda
+\rh_k)}.$$
\endproclaim

Suppose in particular that the $k_\alpha$ are non-negative integers.
Then (12.6) takes the form
$$\vert P_\lambda\vert^2 = 
\prod_{{\alpha\in R^+}}\prodl _{i=0} ^{k_\alpha-1}\frac
{1-q_\alpha^{\langle\lambda + \rh_k,\alpha^\lor \rangle + \frac {1}
{2}k_{\alpha/2} + i}} {1-q_\alpha^{\langle\lambda
+\rh_k,\alpha^\lor \rangle -\frac {1} {2}k_{\alpha/2}-i}}.\tag
$12.6'$ $$

Even when $\lambda = 0$ (so that $P_\lambda =1$) there is something to
be proved here. Indeed, when $R$ is reduced and $\lambda = 0$, the
conjecture ($12.6'$) reduces to the constant term conjectures of
\cite{10} and \cite{13}. For conjecture $A'$ of \cite{13} asserts that
the constant term of the product
$$\prod_{{\alpha\in R^+}}\prod _{i=0}
^{k_\alpha-1}(1-q_\alpha^i e^\alpha)(1-q^{i + 1}_\alpha e^{-\alpha })\tag
12.7$$
should be equal to
$$\prod_{{\alpha\in R}}\frac
{(q_\alpha;q_\alpha)_{\vert\langle\rh_k,\alpha^\lor \rangle+k_\alpha
\vert}}
{(q_\alpha;q_\alpha)_{\vert\langle\rh_k,\alpha^\lor \rangle\vert}}.\tag
12.8$$
Now the product (12.7) is precisely the product $\triangle'$ defined
in (3.10), and by (3.12) the constant term of $\triangle'$ is
$W(t)\vert 1\vert^2$, which by ($12.6'$) is equal to
$$W(t)\prod_{{\alpha\in R^+}}\frac {(q_\alpha;
q_\alpha)_{\langle\rh_k,\alpha^\lor \rangle + k_\alpha -
1}\,(q_\alpha;q_\alpha)_{\langle\rh_k,\alpha^\lor\rangle-k_\alpha}}
{(q_\alpha;q_\alpha)_{\langle\rh_k,\alpha^\lor \rangle
-1}\,(q_\alpha;q_\alpha)_{\langle\rh_k,\alpha^\lor \rangle}}.\tag12.9$$
On the other hand we have
$$W(t)=\prod_{{\alpha\in R^+}}\frac
{1-q_\alpha^{\langle\rh_k,\alpha^\lor \rangle+k_\alpha}}
{1-q_\alpha^{\langle\rh_k,\alpha^\lor \rangle}},$$
by (\cite{8}, 2.4 nr) applied to the root system $S^\lor $. Hence
(12.9) is equal to (12.8), which proves our assertion.

If however $R$ is of type $BC_n$, there are two choices for $S$ when
$n\geq 2$, so that our conjecture ($12.6'$) when $\lambda = 0$ contains
two distinct constant-term conjectures related to the root system
$BC_n$. Neither of these is obviously equivalent to Morris's
Conjecture $A$ for $BC_n$ (\cite{13}, 3.4).

Recall next that each $f\in A$ is regarded as a function on $V$, by
the rule $e^\lambda(x)=q^{\langle\lambda,x\rangle}$ ($x\in V$, 
$\lambda\in P$). Also let
$$\alpha^\ast = (\alpha_\ast)^\lor  = u_\alpha\alpha^\lor $$
for each $\alpha\in R$, and
$$\rh_k^\ast =\frac {1} {2}\sum_{\alpha\in R^+}k_\alpha\alpha^\ast.$$

\proclaim{Conjecture (12.10)} For $\lambda\in P^{++}$,
$$P_\lambda(\rh^\ast_k)=q^{-\langle\lambda,\rh_k^\ast\rangle}c(\rh_k)/
c(\lambda +\rh_k).$$
\endproclaim

Both conjectures (12.6) and (12.10) are true in each of the situations
considered in Sections~8--11:

\medskip\noindent
(i) When $R$ is reduced and the $k_\alpha$ are all equal to $1$,
we have $P_\lambda=\chi_\lambda$ (8.3), and (12.6) reduces to
$\vert\chi_\lambda\vert^2 =1$, i.e., to (8.4). As to (12.10), we have
$\rh^\ast_k = \frac {1} {2}\suml_{\alpha\in R}\alpha^\ast =\rh^\ast$
say, and 
$$P_\lambda (\rh_k^\ast)=\chi_\lambda(\rh^\ast)=\frac
{\suml_{w\in W}\varepsilon(w)\,q^{\langle w(\lambda +
\rh),\rh^\ast\rangle}} {\suml_{w\in W}\varepsilon(w)\,q ^{\langle
w\rh,\rh^\ast\rangle}}.$$
By Weyl's denominator formula for the root system $S^\lor $, this
factorizes to give
$$q^{-\langle\lambda,\rh^\ast\rangle}\prod_{{\alpha\in
R^+}}\frac {1-q_\alpha^{\langle\lambda +
\rh,\alpha^\lor \rangle}} {1-q_\alpha^{\langle\rh,\alpha^\lor \rangle}}$$
in agreement with (12.10).

\smallskip\noindent
(ii) When the $k_\alpha$ are all zero we have $P_\lambda =m_\lambda$
and $\vert P_\lambda\vert^2 = \vert W_\lambda\vert^{-1}$. On the other
hand, it follows from our definitions that when the $k_\alpha$ are all
zero we have $c_\alpha^\ast(-\lambda-\rh_k;q_\alpha)=1$ for all
$\alpha\in R^+$ and $c_\alpha(\lambda +\rh_k; q_\alpha)=1$ for all
$\alpha\in R^+$ such that $\langle\lambda,\alpha^\lor \rangle\ne 0$.
When $\langle\lambda,\alpha^\lor \rangle = 0$ we have to interpret
$c_\alpha(\lambda +\rh_k;q_\alpha)$ by means of the limit relation
(12.1); with $k_\alpha = k$ for all $\alpha$ this leads to
$$\lim_{{k\rightarrow 0}} \ \frac {c^\ast(-\lambda -\rh_k)}
{c(\lambda
+\rh_k)}=\underset{\langle\lambda,\alpha^\lor \rangle=0}\to{\prod_{
{\alpha\in R^+}}}\frac {\langle\rh,\alpha^\lor \rangle + e_\alpha}
{\langle\rh,\alpha^\lor \rangle + e_\alpha + 1}$$
where $e_\alpha = \frac {1} {2}$ if $\frac {1} {2}\alpha\in R$, and
$e_\alpha = 0$ otherwise; and this product is equal to $\vert
W_\lambda\vert^{-1}$. This checks (12.6) in this case, and (12.10) is
analogous (both sides are equal to $\vert W\vert /\vert
W_\lambda\vert$).

\smallskip\noindent
(iii) When $R$ is of rank 1 (\S9) the formulas (9.8), (9.9), (9.14)
and (9.15) show that both conjectures are true.

\smallskip\noindent
(iv) When $R$ is of type $A_n$, the polynomials $P_\lambda$ are
essentially the same as the symmetric functions $P_\lambda(x;q,t)$
studied in \cite{11}, Chapter VI. Both conjectures are true in this
case, and are proved in {\it loc\. cit.}

\smallskip\noindent
(v) In the situation of \S10 we must express everything in terms of
the $t_\alpha$ before setting $q=0$. We have
$$q_\alpha^{\langle\lambda + \rh_k,\alpha^\lor \rangle} =
q_\alpha^{\langle\lambda,\alpha^\lor \rangle}q^{\langle\rh_k,\alpha^\ast
\rangle}$$
and $\langle\rh_k,\alpha^\ast\rangle = \frac {1} {2}\suml_{\beta\in
R^+} k_\beta u_\beta\langle\beta_\ast,\alpha^\ast\rangle$, so that
$$q^{\langle\rh_k,\alpha^\ast\rangle }=\prod_{{\beta\in
R^+}}t_\beta^{\langle\beta_\ast,\alpha^\ast\rangle/2} =
\bold {t}^{\HT(\alpha^\ast)}$$
in the notation of \cite{8}. It follows that
$c^\ast(-\lambda-\rh_k)=1$ when $q=0$, and that
$$c(\lambda+\rh_k)=\underset{\langle\lambda,\alpha^\lor \rangle=0}\to
{\prod_{{\alpha\in R^+}}}
\frac {1-t_{\alpha/2}t_\alpha\bold {t}^{\HT(\alpha^\ast)}}
{1-t_{\alpha/2}\bold {t}^{\HT(\alpha^\ast)}}$$
which from the results of \cite{8} is easily seen to be equal to the
polynomial $W_\lambda(t)$. Hence in the present situation the
right-hand side of (12.6) is equal to $W_\lambda(t)^{-1}$, which by
(10.5) is equal to $\vert P_\lambda\vert^2$.

Next consider (12.10). We have
$$e^\alpha(\rh_k^\ast) =
q^{\langle\alpha,\rh^\ast_k\rangle}=\prod_{{\beta\in
R^+}}t_\beta^{\langle\beta^\lor ,\alpha\rangle/2}=\bold
{t}^{\HT(\beta)}$$
and in the formula (10.1) for $P_\lambda$, when we evaluate at
$\rh_k^\ast$, all the terms will vanish except that corresponding to
$w_0$, the longest element of $W$. Consequently we obtain
$$(e^\lambda P_\lambda)(\rh^\ast_k)=W(t)/W_\lambda(t)$$
which from above is also equal to $c(\rh_k)/c(\lambda +\rh_k)$, thus
verifying (12.10) in this case.

\smallskip\noindent
(vi) Finally, in the ``limiting case" $q\rightarrow 1$ considered in
\S11, Heckman \cite{4} has proved that $\vert
P_\lambda\vert^2_k/\vert 1\vert^2_k$ has the value predicted by (12.6),
and recently Opdam \cite{15} has evaluated $\vert 1\vert^2_k$. These
results confirm that (12.6) is true in the limiting case.

As to (12.10), it follows from the work of Harish-Chandra (\cite{5}
Chapter V) that in the symmetric space situation (11.16) the zonal
spherical function $\omega_\lambda$ for $\lambda\in P^{++}$ is given by
$$\omega_\lambda(\exp  X)=\frac {c(\lambda +\rh_k)}
{c(\rh_k)}P_\lambda(X)\quad \quad \quad (X\in \frak a^+).$$
Hence by comparison with (11.17) we have
$$P_\lambda(0)=c(\rh_k)/c(\lambda +\rh_k)\tag 12.11$$
which proves (12.10) in this case. Recently Opdam  (\cite{15}, 
Cor.~5.2) has proved (12.11) for arbitrary values of the $k_\alpha$.

\enddocument